\documentclass{llncs}

\usepackage{latexsym}
\usepackage{amsmath}
\usepackage{amssymb}
\usepackage{bussproofs}
\usepackage[T1]{fontenc}
\usepackage[latin1]{inputenc}

\newcommand{\NK}{\mathsf{NK}}
\newcommand{\LK}{\mathsf{LK}}
\newcommand{\lm}{{\lambda\mu}}
\newcommand{\mutilde}{{\tilde{\mu}}}
\newcommand{\lmm}{{\bar{\lambda}\mu\mutilde}}
\newcommand{\dagcirc}{{\dag\circ}}
\newcommand{\circdag}{{\circ\dag}}
\newcommand{\trans}[2][\dag]{{#2}^{#1}}
\newcommand{\subst}[3]{{#1 [ #2 \gets #3 ]}}
\newcommand{\reduce}[2]{\ifcase #2\mathop{{\overset{*}{\to}}_{#1}}\else\mathop{\to_{#1}}\fi}
\newcommand{\lreduce}[2]{\ifcase #2\ \mathop{{\overset{*}{\leadsto}}_{#1}}\ \else\ \mathop{\leadsto_{#1}}\ \fi}
\newcommand{\var}[1]{{#1}}
\newcommand{\abs}[2]{{\lambda #1 . #2}}
\newcommand{\app}[2]{{( #1 ) \, #2}}
\newcommand{\mabs}[2]{{\mu #1 . #2}}
\newcommand{\mapp}[2]{{[ #1 ] \, #2}}
\newcommand{\hcst}[1]{{ \{ #1 \}}}
\newcommand{\hhole}[1]{{\bullet}}
\newcommand{\hvar}[1]{{#1}}
\newcommand{\hpush}[2][\beta]{{#1 (#2)}}
\newcommand{\hpop}[2]{\eapp{#2}{#1}}
\newcommand{\happ}[2]{\eapp{#2}{#1}}
\newcommand{\evar}[1]{{#1}}
\newcommand{\eapp}[2]{{#1 \cdot #2}}
\newcommand{\eabs}[2]{{\mutilde #1 . #2}}
\newcommand{\cut}[2]{{\langle #1 \, | \, #2 \rangle}}

\EnableBpAbbreviations
\newcommand{\nullary}{\AXC}
\newcommand{\unary}{\UIC}
\newcommand{\binary}{\BIC}

\newcommand{\rlabel}{\RL}

\newenvironment{prftree}{\everymath{\scriptstyle}$\bgroup}{\DisplayProof\egroup$\everymath{\textstyle}}

\newenvironment{circtree}{\everymath{\scriptstyle}$\bgroup}{\DisplayProof\egroup^\circ$\everymath{\textstyle}}

\title{$\lm$-calculus and duality: \\ call-by-name and call-by-value}

\author{Jérôme Rocheteau}

\institute
{
  \begin{tabular}{c}
    ESTAS -- INRETS \\ 
    20 rue Élisée Reclus -- BP 317 \\
    F-59666 Villeneuve d'Ascq Cedex \\
    jerome.rocheteau@inrets.fr
  \end{tabular}
  \hfill
  \begin{tabular}{c}
    Preuves, Programmes et Systèmes \\
    CNRS -- Université de Paris VII \\
    UMR 7126 -- Case 7014 \\
    175 rue du Chevaleret -- 75013 Paris -- France \\
  \end{tabular}
}

\begin{document}

\maketitle

\begin{abstract}
Under the extension of Curry-Howard's correspondence to classical logic, 
Gentzen's $\NK$ and $\LK$ systems can be seen as syntax-directed systems of simple types 
respectively for Parigot's $\lm$-calculus and Curien-Herbelin's $\lmm$-calculus. 
We aim at showing their computational equivalence. 
We define translations between these calculi. 
We prove simulation theorems for an undirected evaluation as well as for 
call-by-name and call-by-value evaluations. 
\end{abstract}

\section{Introduction}
\label{sec:introduction}

Key systems for classical logic in proof theory are Gentzen's $\NK$ and $\LK$. 
The logical equivalence between the latter was proved in \cite{Gentzen34}. 
We deal with the extension of Curry-Howard's correspondence 
between proofs and programs through the systems of simple types for the $\lm$ and $\lmm$-calculi. 
This extension concerns some other calculi. It is initially Felleisen's $\lambda c$-calculus. 
Its type system is the intuitionistic natural deduction with the double negation axiom. 
Griffin proposed this axiom as the type for the $c$-operator in \cite{Griffin90}. 
However, we focus on calculi that correspond closer to Gentzen's systems. 
The $\lm$-calculus was defined for $\NK$ in \cite{Parigot92}. 
The $\lmm$-calculus was designed for $\LK$ in \cite{CurienHerbelin00}. 
In the general case, these two calculi are not deterministic. There exists critical pairs. 
The $\lmm$-calculus admits two deterministic projections depending on choosing  
one of the two possible symmetric orientations of a critical pair. 
They correspond to the call-by-name/call-by-value duality. 
% This was pointed out by two possible translations of the implication $A \to B$ 
% into classical logic \cite{Wadler03} as well as into polarised linear logic \cite{LaurentRegnier03}. 

% Relations between call-by-name and call-by-value disciplines and $\lambda$-calculi 
% provided with some control operators are based on the notion of evaluation contexts. 
% This helps to formalise an evaluation order \cite{Bierman98}. 
% A context can be seen as a term with a single hole, written $e \hcst{}$. 
% Filling a context hole with a term is written $e \hcst{t}$. 
% A context represents the rest of a computation after the reduction of $t$. 
% In this sense, $e$ is the current continuation of $t$. 

% Firstly, call-by-name and call-by-value simulations 
% into the pure $\almbda$-calculus use continuations \cite{Plotkin75}. 

We aim at proving the computational equivalence between $\lm$ and $\lmm$-calculi. 
A major step was reached with the proof of the simulation 
of the $\lm$-calculus by the $\lmm$-calculus in \cite{CurienHerbelin00}. 
It holds both for call-by-name and call-by-value evaluations. 
We present the call-by-name/call-by-value projections of the $\lm$-calculus 
in the same way as for the $\lmm$ in \cite{CurienHerbelin00} . 
It consists of choosing one of the two possible orientations of a critical pair. 
We prove that the $\lm$-calculus simulates backwards 
the $\lmm$-calculus in such a way that we obtain easily the same result for 
the call-by-name, for the call-by-value and for the simple type case. 
The $\lmm$-calculus is composed of three syntactic categories: 
terms, contexts (or environments) and commands. 
The $\lm$-calculus is basically composed of terms and commands. 
We add contexts to the $\lm$-calculus. 
It eases mappings between the $\lm$ and $\lmm$-calculi. 
We extend the translation from the $\lm$-calculus to 
the $\lmm$-calculus defined in \cite{CurienHerbelin00} over the $\lm$-contexts. 
We define backwards a translation from the $\lmm$-calculus to the $\lm$-calculus. 

In section \ref{sec:lm-calculus} we present the $\lm$-calculus. 
In section \ref{sec:lmm-calculus} we present the $\lmm$-calculus. 
In section \ref{sec:translations-lm-lmm-calculi} we define translations 
between these two calculi. In section \ref{sec:simulations-lm-lmm-calculi} 
we prove simulation theorems that hold for call-by-name and call-by-value. 

\section{$\lm$-calculus}
\label{sec:lm-calculus}

We follow the definition given in \cite{Parigot92}. 
Firstly, we present the grammar of terms and commands. 
Secondly, we present the system of simple types. 
Thirdly, we present generic reductions and their 
call-by-name and call-by-value projections. 
Fourthly, we extend both the grammar and the type system to the contexts. 

Basically, the $\lm$-calculus is composed of terms and commands. 
They are defined by mutual induction: 
\begin{displaymath}
  t \ ::= \ \var{x} \ | \ \abs{x}{t} \ | \ \app{t}{t} \ | \ \mabs{\alpha}{c}
  \qquad
  c \ ::= \ \mapp{\alpha}{t}
\end{displaymath}
Symbols $x$ range over $\lambda$-variables, symbols $\alpha$ range over $\mu$-variables. 
We note $x \in t$ or $\alpha \in t$ the fact that $x$ or $\alpha$ has 
a free occurrence in $t$. Symbols $\lambda$ and $\mu$ are binders. 
Two terms are equal modulo $\alpha$-equivalence. 

The system of simple types for the $\lm$-calculus is based on two kinds of sequents. 
The first $\Gamma \vdash t : T ~|~ \Delta$ concerns the terms and 
the second $c : (\Gamma \vdash \Delta)$ concerns the commands  
in which $T$ is a simple type obtained by the grammar $T \ ::= \ X ~|~ T \to T$, 
$\Gamma$ is a finite domain application from $\lambda$-variables to simple types 
and $\Delta$ is a finite domain application from $\mu$-variables to simple types. 
$\Gamma , \Gamma'$ denotes the union of the applications $\Gamma$ and $\Gamma'$. 
System rules are: 
\begin{center}
  \begin{prftree}
    \nullary{}
    \unary{$x:A \vdash x:A ~|~ $}
  \end{prftree}
  \hfil
  \begin{prftree}
    \nullary{$\Gamma \vdash t:B ~|~ \Delta$}
    \unary{$\Gamma \setminus \{ x : A \} \vdash \abs{x}{t} : A \to B ~|~ \Delta$}
  \end{prftree}
  \hfil
  \begin{prftree}
    \nullary{$\Gamma \vdash u:A \to B ~|~ \Delta$}
    \nullary{$\Gamma' \vdash v:A ~|~ \Delta'$}
    \rlabel{$(\ast)$}
    \binary{$\Gamma , \Gamma' \vdash \app{u}{v} : B ~|~ \Delta, \Delta'$}
  \end{prftree}
\end{center}
\begin{center}
  \begin{prftree}
    \nullary{$\Gamma \vdash t:A ~|~ \Delta$}
    \rlabel{$(\ast)$}
    \unary{$\mapp{\alpha}{t} : (\Gamma \vdash \Delta , \alpha : A  )$}
  \end{prftree}
  \hfil
  \begin{prftree}
    \nullary{$c : (\Gamma\vdash \Delta )$}
    \unary{$\Gamma \vdash \mabs{\alpha}{c} : A ~|~ \Delta \setminus \{ \alpha : A \}$}
  \end{prftree}
\end{center}
The restriction $(\ast)$ requires that $\Gamma$ and $\Gamma'$ match each other 
on the intersection of their domains. This holds for $\Delta$ and $\Delta'$ too. 

The category of contexts is introduced in order to ease 
comparisons with the homonymous category of the $\lmm$-calculus. 
$\lm$-contexts are defined by mutual induction with the terms: 
\begin{displaymath}
  e \ ::= \ \hvar{\alpha} \ | \ \hpush{t} \ | \ \hpop{e}{t} 
\end{displaymath}
We can see contexts as commands with a hole to fill. 
The first construction $\hvar{\alpha}$ expects a term $t$ 
in order to provide the command $\mapp{\alpha}{t}$. 
The second $\hpush{t}$ expects a term $u$ 
in order to provide the command $\mapp{\beta}{\app{t}{u}}$. 
The last $\hpop{h}{t}$ puts the term $t$ on a stack and expects 
another term to fill the hole. 
\begin{definition}
  \label{def:lm-term-context-cut}
  Let $t$ a term and  $e$ a context. 
  The command $e \hcst{t}$ is defined by induction on $e$: 
  \begin{displaymath}
    e \hcst{t} = 
    \begin{cases} 
      \mapp{\alpha}{t} & \mathsf{if \ } e = \hvar{\alpha} \\ 
      \mapp{\beta}{\app{u}{t}} & \mathsf{if \ } e = \hpush{u} \\
      h \hcst{ \app{t}{u} } & \mathsf{if \ } e = \hpop{h}{u}  \\
    \end{cases}  
  \end{displaymath}
\end{definition}
The type system is extended to another kind of sequents 
$\Gamma ~|~ e : T \vdash \Delta$. The typing rules give the context $e$ 
the type of the term $t$ that fills the hole of $e$: 
\begin{center}
  \begin{prftree}
    \nullary{}
    \unary{$~|~ \hvar{\alpha} : A \vdash \alpha : A$}
  \end{prftree}
  \hfil
  \begin{prftree}
    \nullary{$\Gamma \vdash t : (A \to B) ~|~ \Delta$}
    \unary{$\Gamma ~|~ \hpush{t} : A \vdash \Delta , \beta : B$}
  \end{prftree}
  \hfil
  \begin{prftree}
    \nullary{$\Gamma \vdash t : A ~|~ \Delta$}
    \nullary{$\Gamma' ~|~ e : B \vdash \Delta'$}
    \binary{$\Gamma , \Gamma' ~|~ \hpop{e}{t} : (A \to B) \vdash \Delta , \Delta'$}
  \end{prftree}
\end{center}
A sequent calculus like cut-rule can then be derived in this system 
as a term against context application. 
\begin{lemma}
  \label{thm:nk-cut-rule}
  The rule 
  \begin{prftree}
    \nullary{$\Gamma \vdash t : A ~|~ \Delta$}
    \nullary{$\Gamma' ~|~ e : A \vdash \Delta'$}
    \binary{$e \hcst{t} : (\Gamma , \Gamma' \vdash \Delta , \Delta')$}
  \end{prftree}
  holds in $\lm$.  
\end{lemma}

\begin{proof} By induction on $e$. 
  \begin{itemize}
  \item if $e = \hvar{\alpha}$ then $e \hcst{t} = \mapp{\alpha}{t}$ and 
    \begin{prftree}
      \nullary{$\Gamma \vdash t : A ~|~ \Delta$}
      \unary{$\mapp{\alpha}{t} : (\Gamma \vdash \Delta , \alpha : A)$}
      \noLine
      \unary{}
    \end{prftree}
  \item if $e = \hpush{u}$ then $e \hcst{t} = \mapp{\beta}{\app{u}{t}}$ and 
    \begin{center}
      \begin{prftree}
        \nullary{$\Gamma \vdash u : (A \to B) ~|~ \Delta$}
        \nullary{$\Gamma' \vdash t : A ~|~ \Delta'$}
        \binary{$\Gamma, \Gamma' \vdash \app{u}{t} : B ~|~ \Delta , \Delta'$}
        \unary{$\mapp{\beta}{\app{u}{t}} : 
          (\Gamma,\Gamma'\vdash\Delta,\Delta',\beta : B)$}
      \end{prftree}
    \end{center}
  \item if $e = \hpop{h}{u}$ then $e \hcst{t} = h \hcst{\app{t}{u}}$ and 
    \begin{center}
      \begin{prftree}
        \nullary{$\Gamma \vdash t : (A \to B) ~|~ \Delta$}
        \nullary{$\Gamma' \vdash u : A ~|~ \Delta'$}
        \binary{$\Gamma , \Gamma' \vdash \app{t}{u} : B ~|~ \Delta , \Delta'$}
        \nullary{$\Gamma'' ~|~ h : B \vdash \Delta''$}
        \rlabel{ind. \ hyp.}
        \binary{$h\hcst{\app{t}{u}} : 
          (\Gamma,\Gamma',\Gamma''\vdash\Delta,\Delta',\Delta'')$}
      \end{prftree}
    \end{center}
  \end{itemize}
\end{proof}

\begin{definition}
  Let $t$ a term, $e$ a context and $\alpha$ a $\mu$-variable, 
  The term $\subst{t}{\alpha}{e}$ -- the substitution of $\alpha$ by $e$ in $t$ -- 
  is defined by induction on $t$: 
  \begin{displaymath}
    \subst{t}{\alpha}{e} = 
    \begin{cases}
      \var{x} & \mathsf{if \ } t = \var{x} \\
      \abs{x}{\subst{u}{\alpha}{e}} & \mathsf{if \ } t = \abs{x}{u} \\
      \app{\subst{u}{\alpha}{e}}{\subst{v}{\alpha}{e}} & \mathsf{if \ } t = \app{u}{v} \\
      \mabs{\beta}{\subst{c}{\alpha}{e}} & \mathsf{if \ } t = \mabs{\beta}{c} \\
    \end{cases} 
  \end{displaymath}
  \begin{displaymath}
    \subst{c}{\alpha}{e} = 
    \begin{cases}
      e \hcst{\subst{t}{\alpha}{e}} & \mathsf{if \ } c = \mapp{\alpha}{t} \\
      \mapp{\beta}{\subst{t}{\alpha}{e}} & \mathsf{if \ } c = \mapp{\beta}{t} \\
    \end{cases} 
  \end{displaymath}
\end{definition}

The computation notion is based on reductions. 
We remind one-step reduction rules: 
\begin{displaymath}
  \begin{array}{rcl}
    \app{ \abs{x}{u} }{t} & \reduce{\beta}{1} & \subst{u}{x}{t} \\
    \mabs{\delta}{\mapp{\delta}{t}} & \reduce{\theta}{1} & t \ (\mathsf{if \ } \delta \notin t) 
  \end{array}
  \qquad
  \begin{array}{rcl}
    \app{ \mabs{\alpha}{c} }{t} & \reduce{\mu}{1} & \mabs{\alpha}{\subst{c}{\alpha}{\happ{\alpha}{t}}} \\
    \app{t}{ \mabs{\alpha}{c} } & \reduce{\mu'}{1} & \mabs{\alpha}{\subst{c}{\alpha}{\hpush[\alpha]{t}}} \\
    \mapp{\beta}{\mabs{\alpha}{c}} & \reduce{\rho}{1} & \subst{c}{\alpha}{\beta} 
  \end{array}
\end{displaymath}
The reduction $\reduce{\gamma}{0}$ stands for 
the reflexive and transitive closure of $\reduce{\gamma}{1}$ 
and the reduction $\reduce{}{0}$ stands for the union of $\reduce{\gamma}{0}$ 
for $\gamma \in \{ \beta , \mu , \mu' , \rho , \theta \}$. 

Some of these reductions are linear. 
Both of the $\rho$ and $\theta$-reductions are linear because they correspond to the identity in $\NK$. 
The $\beta$-reduction from the term $\app{\abs{x}{t}}{y}$ is linear because 
it consists of replacing a variable by another variable inside a term.
It corresponds to a normalisation against an axiom rule in $\NK$. 
The $\beta$-reduction from the term $\app{\abs{x}{t}}{u}$ 
where $x$ has a single free occurrence in $t$ is linear too 
because it consists either of substituting a single variable occurrence by any term. 
It corresponds either to a normalisation without a proof-tree branch duplication. 

Reductions $\lreduce{\gamma}{1}$, $\lreduce{\gamma}{0}$ and $\lreduce{}{0}$ 
have the same meanings as in the general case. 
The relation $\approx$ is defined as the reflexive, transitive 
and symmetric closure of $\lreduce{}{0}$. 

There exists a critical pair for computation determinism. 
Applicative terms $\app{\abs{x}{t}}{\mabs{\beta}{d}}$ 
and $\app{\mabs{\alpha}{c}}{\mabs{\beta}{d}}$ 
can be $\beta$ or $\mu'$-rewritten in the first case and 
$\mu$ or $\mu'$-rewritten in the second case. 
We can see the call-by-name and call-by-value disciplines as 
restrictions of the generic reductions. 

The call-by-name evaluation consists of allowing every reduction but the $\mu'$-rule. 
The $\beta$-reduction holds in the first case and the $\mu$-reduction in the second. 
Formally the call-by-name reduction is 
$\reduce{n}{0} = \reduce{}{0} \setminus \reduce{\mu'}{0}$. 

The call-by-value evaluation consists of prohibiting 
$\beta$ and $\mu$-reductions in which the argument is a $\mu$-abstraction. 
Formally we define a subset of terms called values by this grammar: 
$v \ ::= \ x \ | \ \abs{x}{t}$. 
$\beta_v$ and $\mu_v$-reductions are defined instead of generic $\beta$ and $\mu$ ones: 
\begin{displaymath}
  \begin{array}{rclcrcl}
    \app{ \abs{x}{u} }{v} & \reduce{\beta_v}{1} & \subst{u}{x}{v} 
    & \qquad & \app{ \mabs{\alpha}{c} }{v} & \reduce{\mu_v}{1} & \mabs{\alpha}{\subst{c}{\alpha}{\happ{\alpha}{v}}} 
  \end{array}
\end{displaymath}
The call-by-value reduction $\reduce{v}{0}$ is the union of $\reduce{\gamma}{0}$ 
for $\gamma \in \{ \beta_v , \mu_v , \mu' , \rho , \theta \}$. 
Critical pairs are then $\mu'$-rewritten. 

There is another way to define call-by-value into the $\lm$-calculus. 
The solution is detailed in \cite{OngStewart97}. 
It consists of restricting the $\mu'$-rule to values 
instead of the $\mu$: 
\begin{displaymath}
    \app{v}{\mabs{\alpha}{c}} 
    \reduce{\mu'_v}{1} 
    \mabs{\alpha}{\subst{c}{\alpha}{\hpush[\alpha]{v}}}
\end{displaymath}
Formally $\reduce{v}{0}$ becomes the union of $\reduce{\gamma}{0}$ 
for $\gamma \in \{ \beta_v , \mu , \mu'_v , \rho , \theta \}$. 
In fact terms $\app{\abs{x}{t}}{\mabs{\alpha}{c}}$ and 
$\app{\mabs{\alpha}{c}}{\mabs{\alpha'}{c'}}$ are respectively 
$\mu'$ and $\mu$-reduced because $\mabs{\alpha}{c}$ is not a value in these cases. 
However, we follow Curien-Herbelin's call-by-value definition. 

We finish this section by a lemma. It is useful for the section \ref{sec:simulations-lm-lmm-calculi} simulation theorems. 
Any command of the form $e \hcst{\mabs{\alpha}{c}}$ is a redex. 
However, some can not be reduced in call-by-name nor in call-by-value. 

\begin{lemma}
  \label{thm:lm-cut-subst}
  \begin{math}
    e \hcst{\mabs{\alpha}{c}} \reduce{}{0} \subst{c}{\alpha}{e}
  \end{math}
\end{lemma}

\begin{proof} By induction on $e$. 
  \begin{itemize}
  \item if $e = \hvar{\beta}$ then $e \hcst{\mabs{\alpha}{c}} = \mapp{\beta}{\mabs{\alpha}{c}} \ \lreduce{\rho}{1} \ \subst{c}{\alpha}{\beta} $
  \item if $e = \hpush{t}$ then 
    \begin{displaymath}
      \begin{array}{rcl}
        e \hcst{\mabs{\alpha}{c}} 
        & = & \mapp{\beta}{\app{t}{\mabs{\alpha}{c}}} \\ 
        & \reduce{\mu'}{1} & \mapp{\beta}{\mabs{\alpha}{\subst{c}{\alpha}{\hpush[\alpha]{t}}}} \\
        & \lreduce{\rho}{1} & \subst{\subst{c}{\alpha}{\hpush[\alpha]{t}}}{\alpha}{\beta} \\
        & = & \subst{c}{\alpha}{\hpush{t}} \\
      \end{array}
    \end{displaymath}
  \item if $e = \happ{h}{t}$ then 
    \begin{displaymath}
      \begin{array}{rcl} 
        e \hcst{\mabs{\alpha}{c}}
        & = & h \hcst{\app{\mabs{\alpha}{c}}{t}} \\ 
        & \reduce{\mu}{1} & h \hcst{\mabs{\alpha}{\subst{c}{\alpha}{\happ{\alpha}{t}}}} \\ 
        & \reduce{}{0} & \subst{\subst{c}{\alpha}{\happ{\alpha}{t}}}{\alpha}{h} \\ 
        & = & \subst{c}{\alpha}{\happ{h}{t}} \\ 
      \end{array}
    \end{displaymath}
  \end{itemize}
\end{proof}
This lemma does not hold in call-by-name for the $\hpush{t}$ induction case 
because no $\mu'$-rule is allowed. It holds in call-by-value if 
$t$ is a value for the $\happ{t}{h}$ induction case. 

\section{$\lmm$-calculus}
\label{sec:lmm-calculus}

The $\lmm$-calculus has the same relation against $\LK$ 
as the $\lm$-calculus against $\NK$. 
Reductions of $\lmm$-calculus correspond to the cut elimination steps in $\LK$ 
as well as the $\lm$-reductions correspond to the $\NK$-normalisation. 
We follow the definition given in \cite{CurienHerbelin00}. 
Firstly, we present the grammar of the $\lmm$-calculus. 
Secondly, we present the simple type system. 
Thirdly, we present generic reductions and their call-by-name and call-by-value projections. 

The $\lmm$-calculus is basically composed of terms, commands and contexts\footnote{In \cite{DoughertyGhilezanLescanneL04} 
  these are referred to respectively callers, callees and capsules. We kept the terminology in \cite{CurienHerbelin00} 
  that sounds closer to its meaning: terms are programs, contexts are environments and commands represent 
  "a closed system containing both the program and its environment".}. 
They are defined by mutual induction: 
\begin{displaymath}
  t \ ::= \ \var{x} \ | \ \abs{x}{t} \ | \ \mabs{\alpha}{c}
  \qquad
  c \ ::= \ \cut{t}{e}
  \qquad
  e \ ::= \ \evar{\alpha} \ | \ \eapp{t}{e} \ | \ \eabs{x}{c}
\end{displaymath}
As in the $\lm$, symbols $x$ range over $\lambda$-variables, 
symbols $\alpha$ range over $\mu$-variables and 
symbols $\lambda$, $\mu$ and $\mutilde$ are binders. 
Terms are equal modulo $\alpha$-equivalence. 

This calculus symmetry looks like $\LK$'s left/right symmetry. 
It is confirmed by its system of simple types. 
This system shares types with the $\lm$-calculus. 
It shares the same kinds of sequents too. Its rules are: 
\begin{center}
  \begin{prftree}
    \nullary{}
    \unary{$ x:A \vdash x:A ~|~ $}
  \end{prftree}
  \hfil
  \begin{prftree}
    \nullary{}
    \unary{$~|~ \alpha : A \vdash \alpha : A $}
  \end{prftree}
  \hfil
  \begin{prftree}
    \nullary{$\Gamma \vdash t:B ~|~ \Delta$}
    \unary{$\Gamma \setminus \{ x : A \} \vdash \abs{x}{t} : A \to B ~|~ \Delta$}
  \end{prftree}
  \hfil
  \begin{prftree}
    \nullary{$\Gamma \vdash t : A ~|~ \Delta$}
    \nullary{$\Gamma' ~|~ e : B \vdash \Delta'$}
    \rlabel{$(\ast)$}
    \binary{$\Gamma , \Gamma' ~|~ \eapp{t}{e} : A \to B \vdash \Delta , \Delta' $}
  \end{prftree}
\end{center}
\begin{center}
  \begin{prftree}
    \nullary{$c : (\Gamma\vdash \Delta )$}
    \unary{$\Gamma \vdash \mabs{\alpha}{c} : A ~|~ \Delta \setminus \{ \alpha:A \} $}
  \end{prftree}
  \hfil
  \begin{prftree}
    \nullary{$c : (\Gamma \vdash \Delta )$}
    \unary{$\Gamma \setminus \{ x:A \} ~|~ \eabs{x}{c} : A \vdash \Delta$}
  \end{prftree}
  \hfil
  \begin{prftree}
    \nullary{$\Gamma \vdash t : A ~|~ \Delta$}
    \nullary{$\Gamma' ~|~ e : A \vdash \Delta'$}
    \rlabel{$(\ast)$}
    \binary{$\cut{t}{e} : (\Gamma , \Gamma' \vdash \Delta , \Delta')$}
  \end{prftree}
\end{center}
The restriction $(\ast)$ is the same as that of $\lm$. 

We present one-step reduction rules. 
Substitutions inside the $\lmm$-calculus are supposed to be known. 
Each rule concerns a command but the $\theta$-rule: 
\begin{displaymath}
  \begin{array}{rclcrcl}
    \cut{\abs{x}{u}}{\eapp{t}{e}} & \reduce{\beta}{1} & \cut{t}{\eabs{x}{\cut{u}{e}}} 
    & \qquad & \cut{\mabs{\alpha}{c}}{e} & \reduce{\mu}{1} & \subst{c}{\alpha}{e} \\
    \mabs{\delta}{\cut{t}{\delta}} & \reduce{\theta}{1} & t \ (\delta \notin t) 
    & \qquad & \cut{t}{\eabs{x}{c}} & \reduce{\mutilde}{1} & \subst{c}{x}{t}
  \end{array}
\end{displaymath}
$\mu$ and $\mutilde$-reductions are duals of each other. 
They correspond to the structural rules in $\LK$. 
Reductions $\reduce{\gamma}{0}$ and $\lreduce{\gamma}{0}$ 
have the same meanings as in the $\lm$-calculus. 
The $\beta$-rule is a mere term modification without term duplication. 
It is therefore a linear reduction. The $\theta$-reduction is linear too. 
There is no $\rho$-reduction. It is a $\mu$-rule particular case 
in which $e = \beta$. 

This system is not deterministic. There is a single critical pair 
$\cut{\mabs{\alpha}{c}}{\eabs{x}{d}}$. It can be both $\mu$ or $\mutilde$-rewritten 
so that Church-Rosser's property does not hold. In fact 
$\cut{\mabs{\alpha}{\cut{x}{\eapp{y}{\alpha}}}}{\eabs{x}{\cut{z}{\eapp{x}{\beta}}}}$ 
is $\mu$-rewritten as $\cut{x}{\eapp{y}{\eabs{x}{\cut{z}{\eapp{x}{\beta}}}}}$ and is 
$\mutilde$-rewritten as 
$\cut{z}{\eapp{\mabs{\alpha}{\cut{x}{\eapp{y}{\alpha}}}}{\beta}}$. 
These are two different normal forms. 

Call-by-name and call-by-value disciplines still deal with this problem. 
They both consist of restricting the context construction. 
The first new grammar is called $\lmm_T$ and the second is called $\lmm_Q$. 

The call-by-name evaluation consists of restricting the $\mu$-rule 
to a subset of contexts that are called stacks. 
$\lmm_T$-grammar is: 
\begin{displaymath}
  t \ ::= \ \var{x} \ | \ \abs{x}{t} \ | \ \mabs{\alpha}{c}
  \quad
  c \ ::= \ \cut{t}{e}
  \quad
  s \ ::= \ \alpha \ | \ \eapp{t}{s}
  \quad
  e \ ::= \ s \ | \ \eabs{x}{c}
\end{displaymath}
The $\mu_n$-rule is restricted to the stacks: 
\begin{displaymath}
  \cut{\mabs{\alpha}{c}}{s} \reduce{\mu_n}{1} \subst{c}{\alpha}{s}
\end{displaymath} 
Call-by-name reduction $\reduce{n}{0}$ is the union of $\reduce{\gamma}{0}$ 
for $\gamma \in \{ \beta , \mu_n , \mutilde , \theta \}$. 
The critical pair can then only be $\mutilde$-rewritten. 
This reduction was proved confluent and stable in the $\lmm_T$-calculus  
in \cite{CurienHerbelin00}. 

The call-by-value oriented grammar consists of allowing 
the $\eapp{t}{e}$ context construction only for values. 
$\lmm_Q$-grammar is: 
\begin{displaymath}
  t \ ::= \ \var{x} \ | \ \abs{x}{t} \ | \ \mabs{\alpha}{c}
  \quad
  v \ ::= \ x \ | \ \abs{x}{t} 
  \quad
  c \ ::= \ \cut{t}{e}
  \quad
  e \ ::= \ \alpha \ | \ \eapp{v}{e} \ | \ \eabs{x}{c}
\end{displaymath}
The $\mutilde_v$-rule is restricted to values: 
\begin{displaymath}
  \cut{v}{\eabs{x}{c}} \reduce{\mutilde_v}{1} \subst{c}{x}{v}
\end{displaymath}
Call-by-value reduction $\reduce{v}{0}$ is the union of $\reduce{\gamma}{0}$ 
for $\gamma \in \{ \beta , \mu , \mutilde_v , \theta \}$. 
The command $\cut{\mabs{\alpha}{c}}{\mabs{\alpha'}{c'}}$ can 
then only be $\mu$-rewritten. 
This reduction was proved confluent and stable in the $\lmm_Q$-calculus  
in \cite{CurienHerbelin00}. 

The $\beta'$-rule contracts as shortcut for both 
a linear $\beta$-rule and a $\mutilde$-rule: 
\begin{displaymath}
  \begin{array}{rcl}
    \cut{\abs{x}{u}}{\eapp{t}{e}} & \reduce{\beta'}{1} & \cut{\subst{u}{x}{t}}{e}
  \end{array}
\end{displaymath}
This $\beta'$-rule is obviously compatible with the call-by-name evaluation. 
It is also compatible with the call-by-value because $t$ is a value 
by definition of $\lmm_Q$.

\section {Translations between $\lm$ and $\lmm$-calculi}
\label{sec:translations-lm-lmm-calculi}

We define a translation $\trans[\dag]{(~)}$ from $\lm$ to $\lmm$. 
It extends that of Curien-Herbelin to the $\lm$-contexts. 
We define backwards a translation $\trans[\circ]{(~)}$ from $\lmm$ to $\lm$. 
We prove properties about their compatibilities with the simple type system 
and about their compositions. 

\begin{definition}
  \label{def:dag-translation}
  Application $\trans[\dag]{(~)}$ maps any $\lm$-term $t$, command $c$ and context $e$ 
  respectively to a $\lmm$-term, command and context.  
  $\trans[\dag]{(~)}$ is defined by induction on $t$, $c$ and $e$: 
  \begin{displaymath}
    \trans[\dag]{t} = 
    \begin{cases}
      \var{x} & \mathsf{ if \ } t = \var{x} \\
      \abs{x}{\trans{u}} & \mathsf{ if \ } t = \abs{x}{u} \\
      \mabs{\beta}{\cut{\trans[\dag]{v}}
        {\eabs{y}{\cut{\trans[\dag]{u}}{\eapp{y}{\beta}}}}} 
      & \mathsf{ if \ } t = \app{u}{v} \ \ (\star) \\
      \mabs{\alpha}{\trans[\dag]{c}} & \mathsf{ if \ } t = \mabs{\alpha}{c} 
    \end{cases}
  \end{displaymath}
  \begin{displaymath}
    \trans{c} = \trans[\dag]{\mapp{\alpha}{t}} = \cut{\trans[\dag]{t}}{\alpha}
  \end{displaymath}
  \begin{displaymath}
    \trans[\dag]{e} = 
    \begin{cases}
      \evar{\alpha} & \mathsf{ if \ } e = \hvar{\alpha}  \\
      \eabs{y}{\cut{\trans[\dag]{t}}{\eapp{y}{\beta}}} & 
      \mathsf{ if \ } e = \hpush{t} \ \ (\star\star) \\
      \eapp{\trans[\dag]{t}}{\trans[\dag]{h}} & \mathsf{ if \ } e = \hpop{h}{t}
    \end{cases}
  \end{displaymath}
\end{definition}
Condition $(\star)$ requires that variables $y$ and $\beta$ have no 
free occurrence in $u$ neither in $v$. 
Condition $(\star\star)$ requires that $y \notin t$. 
A straightforward induction leads us to state that 
$t$ and $\trans[\dag]{t}$ have the same free variables set. 

It seems more natural to translate $\app{u}{v}$ by 
$\mabs{\beta}{\cut{\trans[\dag]{u}}{\eapp{\trans[\dag]{v}}{\beta}}}$.  
This shorter term corresponds in $\LK$ to the arrow elimination rule in $\NK$ too. 
% \begin{flushleft}
%   \begin{dagtree}
%     \nullary{$\Gamma  \vdash u : (A  \to B) ~|~ \Delta $}
%     \nullary{$\Gamma' \vdash v : A  ~|~ \Delta'$}
%     \binary{$\Gamma  , \Gamma' \vdash \app{u}{v} : B ~|~ \Delta  , \Delta'$}
%   \end{dagtree}
%   =
% \end{flushleft}
% \begin{flushright}
%   \begin{prftree}
%     \nullary{$\Gamma  \vdash \trans[\dag]{u} : (A  \to B) ~|~ \Delta  $}
%     \nullary{$\Gamma' \vdash \trans[\dag]{v} : A  ~|~ \Delta'$}
%     \nullary{}
%     \unary{$~|~ \beta : B \vdash \beta : B$}
%     \binary{$\Gamma' ~|~ \eapp{\trans[\dag]{v}}{\beta} 
%       : ( A  \to B ) \vdash \beta : B , \Delta'$}
%     \binary{$\cut{\trans[\dag]{u}}{\eapp{\trans[\dag]{v}}{\beta}} 
%       : (\Gamma  , \Gamma' \vdash \beta 
%       : B , \Delta  , \Delta')$}
%     \unary{$\Gamma  , \Gamma' 
%       \vdash \mabs{\beta}{\cut{\trans[\dag]{u}}{\eapp{\trans[\dag]{v}}{\beta}}}
%       : B ~|~ \Delta  , \Delta$}
%   \end{prftree}
% \end{flushright}
But it would not be compatible with the call-by-value evaluation. 
For example, $\app{x}{\mabs{\alpha}{c}}$ would be translated as 
$\mabs{\beta}{\cut{x}{\eapp{\mabs{\alpha}{\trans[\dag]{c}}}{\beta}}}$ in this case. 
It can not be reduced by any rule in the $\lmm$-calculus. 
However, $\app{x}{\mabs{\alpha}{c}}$ can be $\mu'$-reduced 
in the $\lm$-calculus. 

$\app{u}{v}$ should be translated as 
$\mabs{\beta}{\cut{\trans[\dag]{u}}{\eabs{y}
    {\cut{\trans[\dag]{v}}{\eabs{x}{\cut{y}{\eapp{x}{\beta}}}}}}}$ 
with Ong and Stewart's call-by-value definition in \cite{OngStewart97}. 

We show that translation $\trans[\dag]{(~)}$ is compatible with the type system. 
If a typing environment for a term $t$ exists, it holds for $\trans[\dag]{t}$. 

\begin{lemma}
  \label{thm:dag-type-compatible}
  $\Gamma \vdash t : A ~|~ \Delta 
  \ \Longrightarrow \ 
  \Gamma \vdash \trans[\dag]{t} : A ~|~ \Delta$ 
\end{lemma}

\begin{proof} 
  By a straightforward induction on $t$. 
  We show the less than obvious cases. 
  \begin{itemize}
%   \item if $t = \var{x}$ then $\trans[\dag]{t} = \var{x}$ and 
%     \begin{center}
%       \begin{dagtree}
%         \nullary{}
%         \unary{$x:A \vdash x:A ~|~$}        
%       \end{dagtree}
%       = 
%       \begin{prftree}
%         \nullary{}
%         \unary{$x:A \vdash x:A ~|~$}
%       \end{prftree}
%     \end{center}
%   \item if $t = \abs{x}{u}$ then $\trans[\dag]{t} = $ and 
%     \begin{center}
%       \begin{dagtree}
%         \nullary{$\Gamma \vdash u:B ~|~ \Delta$}
%         \unary{$\Gamma \setminus \{ x : A \} \vdash \abs{x}{u} : A \to B ~|~ \Delta$}
%       \end{dagtree}
%       = 
%       \begin{prftree}
%         \nullary{$\Gamma \vdash \trans[\dag]{u} :B ~|~ \Delta$}
%         \unary{$\Gamma \setminus \{ x:A \} \vdash \abs{x}{\trans[\dag]{u}} : A \to B ~|~ \Delta$}
%       \end{prftree}
%     \end{center}
  \item if $t = \app{u}{v}$ then 
    $\trans[\dag]{t} = 
    \mabs{\beta}{\cut{\trans[\dag]{v}}{\eabs{y}{\cut{\trans[\dag]{u}}{\eapp{y}{\beta}}}}}$
    and 
%     \begin{flushleft}
%       \begin{dagtree}
%         \nullary{$\Gamma \vdash u:A \to B ~|~ \Delta$}
%         \nullary{$\Gamma' \vdash v:A ~|~ \Delta'$}
%         \binary{$\Gamma , \Gamma' \vdash \app{u}{v} : B ~|~ \Delta , \Delta'$}
%       \end{dagtree}
%       = 
%     \end{flushleft}
    \begin{center}
      \begin{prftree}
        \nullary{$\Gamma' \vdash \trans[\dag]{v} : A ~|~ \Delta'$}
        \nullary{$\Gamma \vdash \trans[\dag]{u} : A \to B ~|~ \Delta$}
        \nullary{}
        \unary{$y : A \vdash y : A ~|~$}
        \nullary{}
        \unary{$~|~ \beta : B \vdash \beta : B$}
        \binary{$y : A ~|~ \eapp{y}{\beta} : (A \to B) \vdash \beta : B$}
        \binary{$\cut{\trans[\dag]{u}}{\eapp{y}{\beta}}
          :(\Gamma , y:A \vdash \Delta, \beta:B)$}
        \unary{$\Gamma ~|~ \eabs{y}{\cut{\trans[\dag]{u}}{\eapp{y}{\beta}}}
          \vdash \Delta , \beta : B$}
        \binary{$\cut{\trans[\dag]{v}}{\eabs{y}{\cut{\trans[\dag]{u}}{\eapp{y}{\beta}}}} 
          : (\Gamma , \Gamma' \vdash \Delta , \Delta' , \beta : B)$}
        \unary{$\Gamma , \Gamma' \vdash 
          \mabs{\beta}{\cut{\trans[\dag]{v}}{\eabs{y}{\cut{\trans[\dag]{u}}{\eapp{y}{\beta}}}}} : B 
          |~ \Delta , \Delta'$}
      \end{prftree}
    \end{center}
%   \item if $t = \mabs{\alpha}{c}$ then $\trans[\dag]{t} = \mabs{\alpha}{\trans[\dag]{c}}$ and 
%     \begin{center}
%       \begin{dagtree}
%         \nullary{$c : (\Gamma\vdash \Delta )$}
%         \unary{$\Gamma \vdash \mabs{\alpha}{c}:A ~|~ \Delta \setminus \{ \alpha : A \} $}
%       \end{dagtree}
%       = 
%       \begin{prftree}
%         \nullary{$\trans[\dag]{c} : (\Gamma\vdash \Delta )$}
%         \unary{$\Gamma \vdash \mabs{\alpha}{\trans[\dag]{c}}:A ~|~ \Delta \setminus \{\alpha:A\}$}
%       \end{prftree}
%     \end{center}
%   \item if $c = \mapp{\alpha}{t}$ then $\trans[\dag]{c} = \cut{\trans[\dag]{t}}{\alpha}$ and 
%     \begin{center}
%       \begin{dagtree}
%         \nullary{$\Gamma \vdash t:A ~|~ \Delta$}
%         \unary{$\mapp{\alpha}{t} : (\Gamma\vdash \Delta , \alpha : A  )$}
%       \end{dagtree}
%       = 
%       \begin{prftree}
%           \nullary{$\Gamma \vdash \trans[\dag]{t} : A ~|~ \Delta$}
%           \nullary{}
%           \unary{$~|~ \alpha : A \vdash \alpha : A$}
%           \binary{$\cut{\trans[\dag]{t}}{\alpha} : (\Gamma \vdash \Delta , \alpha : A  )$}
%       \end{prftree}
%     \end{center}
%   \item if $e = \hvar{\alpha}$ then $\trans[\dag]{e} = \evar{\alpha}$ and 
%     \begin{center}
%       \begin{dagtree}
%         \nullary{}
%         \unary{$~|~ \hvar{\alpha} : A \vdash \alpha : A$}
%       \end{dagtree}
%       = 
%       \begin{prftree}
%         \nullary{}
%         \unary{$~|~ \evar{\alpha} : A \vdash \alpha : A$}
%       \end{prftree}
%     \end{center}
  \item if $e = \hpush{t}$ then $\trans[\dag]{e} = \eabs{y}{\cut{\trans[\dag]{t}}{\eapp{y}{\beta}}}$ and 
    \begin{center}
%       \begin{dagtree}
%           \nullary{$\Gamma \vdash t : (A \to B) ~|~ \Delta$}
%           \unary{$\Gamma ~|~ \hpush{t} : A \vdash \Delta , \beta : B$}
%       \end{dagtree}
%       = 
      \begin{prftree}
        \nullary{$\Gamma \vdash \trans[\dag]{t} : (A \to B) ~|~ \Delta$}
        \nullary{}
        \unary{$y : A \vdash y : A ~|~$}
        \nullary{}
        \unary{$~|~ \beta : B \vdash \beta : B$}
        \binary{$y : A ~|~ \eapp{y}{\beta} : (A \to B)\vdash \beta : B$}
        \binary{$\cut{\trans[\dag]{t}}{\eapp{y}{\beta}} 
          : (\Gamma , y: A \vdash \Delta , \beta : B)$}
        \unary{$\Gamma ~|~ \eabs{y}{\cut{\trans[\dag]{t}}{\eapp{y}{\beta}}} 
          : A \vdash \Delta , \beta : B$}
      \end{prftree}
    \end{center}
%   \item if $e = \hpop{h}{t}$ then $\trans[\dag]{e} = \eapp{\trans[\dag]{t}}{\trans[\dag]{h}}$ and 
%     \begin{center}
%       \begin{dagtree}
%         \nullary{$\Gamma \vdash t : A ~|~ \Delta$}
%         \nullary{$\Gamma' ~|~ h : B \vdash \Delta'$}
%         \binary{$\Gamma , \Gamma' ~|~ \hpop{h}{t} : (A \to B) \vdash \Delta , \Delta'$}
%       \end{dagtree}
%       = 
%       \begin{prftree}
%         \nullary{$\Gamma \vdash \trans[\dag]{t} : A ~|~ \Delta$}
%         \nullary{$\Gamma' ~|~ \trans[\dag]{h} : B \vdash \Delta'$}
%         \binary{$\Gamma , \Gamma' ~|~ \eapp{\trans[\dag]{t}}{\trans[\dag]{h}} : (A \to B) 
%             \vdash \Delta , \Delta'$}
%       \end{prftree}
%     \end{center}
  \end{itemize}
\end{proof}

\begin{definition}
  \label{def:circ-translation}
  Application $\trans[\circ]{(~)}$ maps backwards 
  any $\lmm$-term $t$ to a $\lm$-term.
  Definition \ref{def:lm-term-context-cut} is used 
  to translate any $\lmm$-command $c$. 
  Definition of the $\lm$-contexts is used 
  to map the $\lmm$-contexts $e$ as well. 
  $\trans[\circ]{(~)}$ is built by induction on $t$, $c$ and $e$: 
  \begin{displaymath}
    \trans[\circ]{t} = 
    \begin{cases}
      \var{x} & \mathsf{ if \ } t = \var{x} \\
      \abs{x}{\trans[\circ]{u}} &  \mathsf{ if \ } t = \abs{x}{u} \\
      \mabs{\alpha}{\trans[\circ]{c}} & \mathsf{ if \ } t = \mabs{\alpha}{c} 
    \end{cases}
  \end{displaymath}
  \begin{displaymath}
    \trans[\circ]{c} = \trans[\circ]{\cut{t}{e}} = \trans[\circ]{e}\hcst{\trans[\circ]{t}}
  \end{displaymath}
  \begin{displaymath}
    \trans[\circ]{e} = 
    \begin{cases}
      \hvar{\alpha}  & \mathsf{ if \ } e = \evar{\alpha} \\ 
      \hpop{\trans[\circ]{h}}{\trans[\circ]{t}} & \mathsf{ if \ } e = \eapp{t}{h} \\ 
      \hpush{\abs{x}{\mabs{\delta}{\trans[\circ]{c}}}} & \mathsf{ if \ } e = \eabs{x}{c} 
      \ \ (\ast)
    \end{cases}
  \end{displaymath}
\end{definition}
Condition $(\ast)$ requires that $\delta \notin c$. 
$t$ and $\trans[\circ]{t}$ have the same free variables set. 
Application $\trans[\circ]{(~)}$ is compatible with the type system too. 

\begin{lemma}
  \label{thm:circ-type-compatible}
  $\Gamma \vdash t : A ~|~ \Delta 
  \ \Longrightarrow  \ 
  \Gamma \vdash \trans[\circ]{t} : A ~|~ \Delta$ 
\end{lemma}

\begin{proof} By a straightforward induction on $t$. We give two cases. 
  \begin{itemize}
%   \item if $t = \var{x}$ then $\trans[\circ]{t} = \var{x}$ and 
%     \begin{center}
%       \begin{circtree}
%         \nullary{}
%         \unary{$ x:A \vdash x:A ~|~ $}
%       \end{circtree}
%       =
%       \begin{prftree}
%         \nullary{}
%         \unary{$ x:A \vdash x:A ~|~ $}
%       \end{prftree}
%     \end{center}
%   \item if $t = \abs{x}{u}$ then $\trans[\circ]{t} = \abs{x}{\trans[\circ]{u}}$ and 
%     \begin{center}
%       \begin{circtree}
%         \nullary{$\Gamma \vdash t:B ~|~ \Delta$}
%         \unary{$\Gamma \setminus \{ x : A \} \vdash \abs{x}{t} : A \to B ~|~ \Delta$}
%       \end{circtree}
%       =
%       \begin{prftree}
%           \nullary{$\Gamma \vdash \trans[\circ]{t} : B ~|~ \Delta$}
%           \unary{$\Gamma \setminus \{ x : A \} \vdash \abs{x}{\trans[\circ]{t}} 
%             : A \to B ~|~ \Delta$}
%       \end{prftree}
%     \end{center}
%   \item if $t = \mabs{x}{c}$ then $\trans[\circ]{t} = \mabs{x}{\trans[\circ]{c}}$ and 
%     \begin{center}
%       \begin{circtree}
%         \nullary{$c : (\Gamma\vdash \Delta )$}
%         \unary{$\Gamma \vdash \mabs{\alpha}{c} : A ~|~ \Delta \setminus \{ \alpha:A \} $}
%       \end{circtree}
%       =
%       \begin{prftree}
%         \nullary{$\trans[\circ]{c} : (\Gamma\vdash \Delta )$}
%         \unary{$\Gamma \vdash \mabs{\alpha}{\trans[\circ]{c}} 
%           : A ~|~ \Delta \setminus \{ \alpha:A \} $}
%       \end{prftree}
%     \end{center}
  \item if $c = \cut{t}{e}$ then 
    $\trans[\circ]{c} = \trans[\circ]{e} \hcst{\trans[\circ]{t}}$ and 
    \begin{center}
      \begin{circtree}
        \nullary{$\Gamma \vdash t : A ~|~ \Delta$}
        \nullary{$\Gamma' ~|~ e : A \vdash \Delta'$}
        \binary{$\cut{t}{e} : (\Gamma , \Gamma' \vdash \Delta , \Delta')$}
      \end{circtree}
      =
      \begin{prftree}
        \nullary{$\Gamma \vdash \trans[\circ]{t} : A ~|~ \Delta$}
        \nullary{$\Gamma' ~|~ \trans[\circ]{e} : A \vdash \Delta'$}
        \rlabel{lem. \ref{thm:nk-cut-rule}}
        \binary{$\trans[\circ]{e} \hcst{\trans[\circ]{t}} 
          : (\Gamma , \Gamma' \vdash \Delta , \Delta')$}
      \end{prftree}
    \end{center}
%   \item if $e = \evar{\alpha}$ then $\trans[\circ]{e} = \hvar{\alpha}$ and 
%     \begin{center}
%       \begin{circtree}
%         \nullary{}
%         \unary{$~|~ \evar{\alpha} : A \vdash \alpha : A $}
%       \end{circtree}
%       =
%       \begin{prftree}
%         \nullary{}
%         \unary{$~|~ \hvar{\alpha} : A \vdash \alpha : A $}
%       \end{prftree}
%     \end{center}
%   \item if $e = \eapp{t}{h}$ then $\trans[\circ]{e} = \hpop{\trans[\circ]{h}}{\trans[\circ]{t}}$ and 
%     \begin{center}
%       \begin{circtree}
%         \nullary{$\Gamma \vdash t : A ~|~ \Delta$}
%         \nullary{$\Gamma' ~|~ e : B \vdash \Delta'$}
%         \binary{$\Gamma , \Gamma' ~|~ \eapp{t}{e} : A \to B \vdash \Delta , \Delta' $}
%       \end{circtree}
%       =
%       \begin{prftree}
%         \nullary{$\Gamma \vdash \trans[\circ]{t} : A ~|~ \Delta$}
%         \nullary{$\Gamma' ~|~ \trans[\circ]{e} : B \vdash \Delta'$}
%         \binary{$\Gamma , \Gamma' ~|~ \hpop{\trans[\circ]{h}}{\trans[\circ]{t}} 
%           : A \to B \vdash \Delta , \Delta' $}
%       \end{prftree}
%     \end{center}
  \item if $e = \eabs{x}{c}$ then 
    $\trans[\circ]{e}=\hpush{\abs{x}{\mabs{\delta}{\trans[\circ]{c}}}}$ and 
    \begin{center}
      \begin{circtree}
          \nullary{$c : (\Gamma \vdash \Delta )$}
          \unary{$\Gamma \setminus \{ x:A \} ~|~ \eabs{x}{c} : A \vdash \Delta$}
      \end{circtree}
      =
      \begin{prftree}
        \nullary{$\trans[\circ]{c} : (\Gamma \vdash \Delta )$}
        \unary{$\Gamma \vdash \mabs{\delta}{\trans[\circ]{c}} : B ~|~ \Delta$}
        \unary{$\Gamma \setminus \{ x:A \} \vdash \abs{x}{\mabs{\delta}{\trans[\circ]{c}}} 
          : (A \to B)$}
        \unary{$\Gamma \setminus \{ x:A \} ~|~ 
          \hpush{\abs{x}{\mabs{\delta}{\trans[\circ]{c}}}}
          : A \vdash \Delta , \beta : B$}
      \end{prftree}
    \end{center}
  \end{itemize}
\end{proof}

We focus on properties about the composition of $\trans[\dag]{(~)}$ and $\trans[\circ]{(~)}$. 
We want to state that $\trans[\dagcirc]{t} = t$ and that $\trans[\circdag]{t} = t$ 
for any term. But it is not the case, these results hold modulo linear reductions. 

\begin{theorem}
  \label{thm:lm-lmm-composition}
  $\trans[\dagcirc]{t} \lreduce{}{0} t$ 
\end{theorem}

\begin{proof} 
  By a straightforward induction on $t$. 
  Every cases is obtained successively by 
  expanding definitions 
  \ref{def:dag-translation}, 
  \ref{def:lm-term-context-cut}, 
  \ref{def:circ-translation} 
  and by applying the induction hypothesis. 
  We give the case which uses linear reductions additionally. 
  \begin{itemize}
%   \item if $t = \var{x}$ then $x^\dagcirc = x^\circ = x$ and $x^\dagcirc \lreduce{}{0} x$
%   \item if $t = \abs{x}{u}$ then 
%     $(\abs{x}{u})^\dagcirc = (\abs{x}{\trans[\dag]{u}})^\circ = \abs{x}{\trans[\dag]{u}circ}$, 
%     $\trans[\dag]{u}circ \lreduce{}{0} u$ and $\abs{x}{\trans[\dag]{u}circ} \lreduce{}{0} \abs{x}{u}$
  \item if $t = \app{u}{v}$ then 
    \begin{displaymath}
      \begin{array}{rcl}
        \trans[\dagcirc]{\app{u}{v}} 
        % \mathsf{def. \ \ref{def:dag-translation}} 
        & = & 
        \trans[\circ]{\mabs{\beta}{\cut{\trans[\dag]{v}}
            {\eabs{y}{\cut{\trans[\dag]{u}}{\eapp{y}{\beta}}}}}} \\
        % \mathsf{def. \ \ref{def:circ-translation}} 
        % \mathsf{def. \ \ref{def:lm-term-context-cut}} 
        & = & 
        \mabs{\beta}{\mapp{\gamma}
          {\app{\abs{y}{\mabs{\delta}{\mapp{\beta}
                  {\app{\trans[\dagcirc]{u}}{y}}}}}{\trans[\dagcirc]{v}}}} \\ 
        % \mathsf{ind. \ hyp.} 
        & \lreduce{}{0} & 
        \mabs{\beta}{\mapp{\gamma}
          {\app{\abs{y}{\mabs{\delta}{\mapp{\beta}{\app{u}{y}}}}}{v}}} \\
        & \lreduce{\beta}{1} & 
        \mabs{\beta}{\mapp{\gamma}{\mabs{\delta}{\mapp{\beta}{\app{u}{v}}}}} \\
        & \lreduce{\rho}{1} & 
        \mabs{\beta}{\mapp{\beta}{\app{u}{v}}} \\ 
        & \lreduce{\theta}{1} & \app{u}{v} 
      \end{array}
    \end{displaymath}
%   \item if $t = \mabs{\alpha}{c}$ then 
%     $(\mabs{\alpha}{c})^\dagcirc=(\mabs{\alpha}{\trans[\dag]{c}})^\circ=\mabs{\alpha}{\trans[\dag]{c}circ}$, 
%     $\trans[\dag]{c}circ \lreduce{}{0} c$ and 
%     $\mabs{\alpha}{\trans[\dag]{c}circ} \lreduce{}{0} \mabs{\alpha}{c}$ 
%   \item if $c = \mapp{\alpha}{t}$ then 
%     $(\mapp{\alpha}{t})^\dagcirc = (\cut{\trans[\dag]{t}}{\alpha})^\circ = [\alpha]\trans[\dag]{t}circ$, 
%     $\trans[\dag]{t}circ \lreduce{}{0} t$ and $(\mapp{\alpha}{t})^\dagcirc \lreduce{}{0} [\alpha]t$
  \end{itemize}
\end{proof}

We prove two lemmas before stating backwards that $\trans[\circdag]{(~)}$ 
is the identity modulo linear reductions. 
The first lemma is useful to prove the second. 

\begin{lemma}
  \label{thm:dag-applicative-sequence}
  \begin{math}
    \cut{\trans[\dag]{t_0 t_1 \ldots t_n}}{e} 
    \lreduce{}{0} 
    \cut{\trans[\dag]{t_0}}{\trans[\dag]{t_1} \cdot \ldots \cdot \trans[\dag]{t_n} \cdot e}
  \end{math}
\end{lemma}

\begin{proof} By induction on $n$. 
  \begin{itemize}
  \item if $n = 0$ then it is obvious 
  \item if $n = m + 1$ then 
    \begin{displaymath}
      \begin{array}{rcl}
        \cut{\trans[\dag]{t_0 t_1 \ldots t_m t_{m+1}}}{e} 
        & = & \cut{\mabs{\beta}{\cut{\trans[\dag]{t_{m+1}}}{\eabs{y}{\cut{\trans[\dag]{t_0 t_1 \ldots t_m}}{\eapp{y}{\beta}}}}}}{e} \\
        & \lreduce{\mu}{1} & \cut{\trans[\dag]{t_{m+1}}}{\eabs{y}{\cut{\trans[\dag]{t_0 t_1 \ldots t_m}}{\eapp{y}{e}}}} \\
        & \lreduce{\mutilde}{1} & \cut{\trans[\dag]{t_0 t_1 \ldots t_m}}{\eapp{\trans[\dag]{t_{m+1}}}{e}} \\
        % \mathsf{ind. \ hyp.} 
        & \lreduce{}{0} & \cut{\trans[\dag]{t_0}}{\trans[\dag]{t_1} \cdot \ldots \cdot \trans[\dag]{t_m} \cdot \trans[\dag]{t_{m+1}} \cdot e}
      \end{array}
    \end{displaymath}
  \end{itemize}
\end{proof}

The second lemma shows how to map a definition \ref{def:lm-term-context-cut} command. 

\begin{lemma}
  \label{thm:dag-cut-term}
  \begin{math}
    \trans[\dag]{e \hcst{t}} \lreduce{}{0} \cut{\trans[\dag]{t}}{\trans[\dag]{e}}
  \end{math}
\end{lemma}

\begin{proof} By induction on $e$. 
  \begin{itemize}
  \item if $e = \hvar{\alpha}$ then it is obvious by definitions 
    \ref{def:lm-term-context-cut} and \ref{def:dag-translation} 
%     \begin{displaymath}
%       \begin{array}{rcl}
%         & & \trans[\dag]{\hvar{\alpha} \hcst{t}} \\ 
%         \mathsf{def. \ \ref{def:lm-term-context-cut}} & = & 
%         \trans[\dag]{\mapp{\alpha}{t}} \\
%         \mathsf{def. \ \ref{def:dag-translation}} & = & 
%         \cut{\trans[\dag]{t}}{\alpha} \\
%       \end{array}
%     \end{displaymath}
  \item if $e = \hpush{u}$ then 
    \begin{displaymath}
      \begin{array}{rcl}
        \trans[\dag]{\hpush{u} \hcst{t}} 
        % \mathsf{def. \ \ref{def:lm-term-context-cut}} 
        & = & 
        \trans[\dag]{\mapp{\beta}{\app{u}{t}}} \\
        % \mathsf{def. \ \ref{def:dag-translation}} 
        & = & \cut{\mabs{\gamma}{\cut{\trans[\dag]{t}}{\eabs{y}{\cut{\trans[\dag]{u}}{\eapp{y}{\gamma}}}}}}{\beta} \\
        & \lreduce{\mu}{1} & \cut{\trans[\dag]{t}}{\eabs{y}{\cut{\trans[\dag]{u}}{\eapp{y}{\beta}}}} \\
        % \mathsf{def. \ \ref{def:dag-translation}} 
        & = & \cut{\trans[\dag]{t}}{\trans[\dag]{\hpush{u}}} 
      \end{array}
    \end{displaymath}
  \item if $e = \hpop{h}{u}$ then 
    \begin{displaymath}
      \begin{array}{rcl}
        \trans[\dag]{ \hpop{h}{u} \hcst{t}} 
        % \mathsf{def. \ \ref{def:lm-term-context-cut}} 
        & = & \trans[\dag]{h \hcst{\app{t}{u}}} \\ 
        % \mathsf{ind. \ hyp.} 
        & \lreduce{}{0} & \cut{\trans[\dag]{\app{t}{u}}}{\trans[\dag]{h}} \\
        % \mathsf{lem. \ \ref{thm:dag-applicative-sequence}} 
        & \lreduce{}{0} & \cut{\trans[\dag]{t}}{\eapp{\trans[\dag]{u}}{\trans[\dag]{h}}} \\
        % \mathsf{def. \ \ref{def:dag-translation}} 
        & = & \cut{\trans[\dag]{t}}{\trans[\dag]{ \hpop{h}{u}}} 
      \end{array}
    \end{displaymath}
  \end{itemize}
\end{proof}

\begin{theorem}
  \label{thm:lmm-lm-composition}
  $\trans[\circdag]{t} \lreduce{}{0} t$ 
\end{theorem}

\begin{proof} By induction on $t$. 
  We apply definitions \ref{def:dag-translation}, \ref{def:circ-translation} 
  successively and the induction hypothesis. 
  We give a typical case and another which needs 
  either the previous lemma or linear reductions. 
  \begin{itemize}
%   \item if $t = x$ then $x^\circdag = x^\dag = x$ and $x^\circdag \lreduce{}{0} x$ 
%   \item if $t = \abs{x}{u}$ then 
%     $(\abs{x}{u})^\circdag = (\abs{x}{\trans[\circ]{u}})^\dag = \abs{x}{\trans[\circ]{u}dag}$, 
%     $\trans[\circ]{u}dag \lreduce{}{0} u$ and $\abs{x}{\trans[\circ]{u}dag} \leadsto \abs{x}{u}$
%   \item if $t = \mabs{\alpha}{c}$ then 
%     $(\mabs{\alpha}{c})^\circdag=(\mabs{\alpha}{\trans[\circ]{c}})^\dag=\mabs{\alpha}{\trans[\circ]{c}dag}$, 
%     $\trans[\circ]{c}dag \lreduce{}{0} c$ and $\mabs{\alpha}{\trans[\circ]{c}dag} \leadsto \mabs{\alpha}{c}$
  \item if $c = \cut{t}{e}$ then 
    \begin{displaymath}
      \begin{array}{rcl}
        \trans[\circdag]{\cut{t}{e}} 
        % \mathsf{def. \ \ref{def:circ-translation}} 
        & = & 
        \trans[\dag]{\trans[\circ]{e} \hcst{\trans[\circ]{t}}} \\ 
        % \mathsf{lem. \ \ref{thm:dag-cut-term}} 
        & \lreduce{}{0} & 
        \cut{\trans[\circdag]{t}}{\trans[\circdag]{e}} \\
        % \mathsf{ind. \ hyp.} 
        & \lreduce{}{0} & \cut{t}{e} 
      \end{array}
    \end{displaymath}
%   \item if $e =\alpha$ then $\alpha^\circdag = \hvar{\alpha}^\dag = \alpha$ 
%     and $\alpha^\circdag \lreduce{}{0} \alpha$
%   \item if $e = \eapp{t}{h}$ then 
%     $(\eapp{t}{h})^\circdag = (\hpop{\trans[\circ]{h}}{\trans[\circ]{t}})^\dag = \eapp{\trans[\circ]{t}dag}{\trans[\circ]{h}dag}$ 
%       and $\eapp{\trans[\circ]{t}dag}{\trans[\circ]{h}dag} \leadsto \eapp{t}{h}$
  \item if $e = \eabs{x}{c}$ then 
    \begin{displaymath}
      \begin{array}{rcl}
        \trans[\circdag]{\eabs{x}{c}} 
        % \mathsf{def. \ \ref{def:circ-translation}} 
        & = & 
        \trans[\dag]{\hpush{\abs{x}{\mabs{\beta}{\trans[\circ]{c}}}}} \\ 
        % \mathsf{def. \ \ref{def:dag-translation}} 
        & = & 
        \eabs{y}{\cut{\abs{x}{\mabs{\beta}{\trans[\circdag]{c}}}}{\eapp{y}{\beta}}} \\
        % \mathsf{ind. \ hyp.} 
        & \lreduce{}{0} & 
        \eabs{y}{\cut{\abs{x}{\mabs{\beta}{c}}}{\eapp{y}{\beta}}} \\ 
        & \lreduce{\beta}{1} & \eabs{y}{\cut{y}{\eabs{x}{\cut{\mabs{\beta}{c}}{\beta}}}}\\
        & \lreduce{\mutilde}{1} & \eabs{x}{\cut{\mabs{\beta}{c}}{\beta}} \\ 
        & \lreduce{\mu}{1} & \eabs{x}{c} 
      \end{array}
    \end{displaymath}
  \end{itemize}
\end{proof}

\section{Simulations between $\lm$ and $\lmm$-calculi}
\label{sec:simulations-lm-lmm-calculi}

We want to prove that the $\lm$-calculus simulates and is simulated backwards by 
the $\lmm$-calculus. We focus on the undirected evaluation. 
Call-by-name and call-by-value are drawn from this. 

We begin with the simulation of the $\lm$ by the $\lmm$. 
The next four lemmas show results of a $\lm$-substitution 
after a $\beta$, $\mu$, $\mu'$ and $\rho$-reduction. 
Each proof consists successively of 
\begin{itemize}
\item expanding the $\lm$-substitution 
\item expanding the definition of $\trans[\dag]{(~)}$
\item applying the induction hypothesis if necessary
\item factorising the $\lmm$-substitution
\item factorising the definition of $\trans[\dag]{(~)}$ 
\end{itemize}
We give basic cases and those which use lemmas additionally for any proof. 

\begin{lemma}
  \label{thm:dag-beta-subst}
  $\trans[\dag]{\subst{t}{x}{u}} = \subst{\trans[\dag]{t}}{x}{\trans[\dag]{u}}$ 
\end{lemma}

\begin{proof} By induction on $t$.
  \begin{itemize}
  \item if $t = x$ then 
    \begin{math}
      \trans[\dag]{\subst{x}{x}{u}} 
      = \trans[\dag]{u} 
      = \subst{\trans[\dag]{x}}{x}{\trans[\dag]{u}}
    \end{math}
%     \begin{displaymath}
%       \begin{array}{rcl}
%         & & \trans[\dag]{\subst{x}{x}{u}} \\
%         \mathsf{\lm-subst. \ def.} & = & \trans[\dag]{u} \\
%         \mathsf{\lmm-subst. \ def.} & = &  \subst{\trans[\dag]{x}}{x}{\trans[\dag]{u}} 
%       \end{array}
%     \end{displaymath}
  \item if $t = y$ then 
    \begin{math}
      \trans[\dag]{\subst{y}{x}{u}} 
      = y 
      = \subst{\trans[\dag]{y}}{x}{\trans[\dag]{u}}
    \end{math}
%     \begin{displaymath}
%       \begin{array}{rcl}
%         & & \trans[\dag]{\subst{y}{x}{u}}  \\
%         \mathsf{\lm-subst. \ def.} & = & \trans[\dag]{y} \\
%         \mathsf{\lmm-subst. \ def.} & = & \subst{\trans[\dag]{y}}{x}{\trans[\dag]{u}}
%     \end{array}
%     \end{displaymath}
%   \item if $t = \abs{y}{v}$ then 
%     $\trans[\dag]{t} = \abs{y}{\trans[\dag]{v}}$, 
%     $\subst{\trans[\dag]{t}}{x}{\trans[\dag]{u}} = \abs{y}{\subst{\trans[\dag]{v}}{x}{\trans[\dag]{u}}}$ and 
%     \begin{displaymath}
%       \begin{array}{rcl}
%         \subst{t}{x}{u}^\dag 
%         & = & \subst{\abs{y}{v}}{x}{u}^\dag \\
%         & = & \abs{y}{\subst{v}{x}{u}}^\dag \\
%         & = & \abs{y}{\subst{\trans[\dag]{v}}{x}{\trans[\dag]{u}}} \\
%       \end{array}
%     \end{displaymath}
  \item if $t = \app{v}{w}$ then 
    \begin{displaymath}
      \begin{array}{rcl}
        \trans[\dag]{\subst{\app{v}{w}}{x}{u}} 
        % \mathsf{\lm-subst. \ def.} 
        & = & 
        \trans[\dag]{\app{\subst{v}{x}{u}}{\subst{w}{x}{u}}} \\ 
        % \mathsf{def. \ \ref{def:dag-translation}} 
        & = & 
        \mabs{\beta}{\cut{\trans[\dag]{\subst{w}{x}{u}}}
          {\eabs{y}{\cut{\trans[\dag]{\subst{v}{x}{u}}}{\eapp{y}{\beta}}}}} \\ 
        % \mathsf{ind. \ hyp.} 
        & = & 
        \mabs{\beta}{\cut{\subst{\trans[\dag]{w}}{x}{\trans[\dag]{u}}}
          {\eabs{y}{\cut{\subst{\trans[\dag]{v}}{x}
                {\trans[\dag]{u}}}{\eapp{y}{\beta}}}}} \\ 
        % \mathsf{\lmm-subst. \ def.} 
        & = & 
        \subst{\mabs{\beta}{\cut{\trans{w}}{\eabs{y}{\cut{\trans{v}}{\eapp{y}{\beta}}}}}} 
        {x}{\trans[\dag]{u}} \\ 
        % \mathsf{def. \ \ref{def:dag-translation}} 
        & = & 
        \subst{\trans[\dag]{\app{v}{w}}}{x}{\trans[\dag]{u}}
      \end{array}
    \end{displaymath}
%   \item if $t = \mabs{\alpha}{c}$ then 
%     $\trans[\dag]{t} = \mabs{\alpha}{\trans[\dag]{c}}$, 
%     $\subst{\trans[\dag]{t}}{x}{\trans[\dag]{u}} = \mabs{\alpha}{\subst{\trans[\dag]{c}}{x}{\trans[\dag]{u}}}$ and 
%     \begin{displaymath}
%       \begin{array}{rcl}
%         \subst{t}{x}{u}^\dag 
%         & = & \subst{\mabs{\alpha}{c}}{x}{u}^\dag \\
%         & = & \mabs{\alpha}{\subst{c}{x}{u}}^\dag \\
%         & = & \mabs{\alpha}{\subst{\trans[\dag]{c}}{x}{\trans[\dag]{u}}} \\
%       \end{array}
%     \end{displaymath}
%   \item if $c = \mapp{\alpha}{t}$ then 
%     $\trans[\dag]{c} = \cut{\trans[\dag]{t}}{\alpha}$, 
%     $\subst{\trans[\dag]{c}}{x}{\trans[\dag]{u}} = \cut{\subst{\trans[\dag]{t}}{x}{\trans[\dag]{u}}}{\alpha}$ and 
%     \begin{displaymath}
%       \begin{array}{rcl}
%         \subst{c}{x}{u}^\dag 
%         & = & \subst{\mapp{\alpha}{t}}{x}{u}^\dag \\
%         & = & \mapp{\alpha}{\subst{t}{x}{u}}^\dag \\ 
%         & = & \cut{\subst{t}{x}{u}^\dag}{\alpha} \\
%         & = & \cut{\subst{\trans[\dag]{t}}{x}{\trans[\dag]{u}}}{\alpha} \\
%       \end{array}
%     \end{displaymath}
  \end{itemize}
\end{proof}

\newcommand{\substA}[4]{\subst{#1}{#2}{\happ{#2}{#4}}}
\newcommand{\evalA}[3]{\subst{#1}{#2}{\eapp{#3}{#2}}}

\begin{lemma}
  \label{thm:dag-mu-subst}
  \begin{math}
    \trans[\dag]{\substA{t}{\alpha}{v}{u}}
    \lreduce{}{0} 
    \evalA{\trans[\dag]{t}}{\alpha}{\trans[\dag]{u}}
  \end{math}
\end{lemma}

\begin{proof} By induction on $t$.
  \begin{itemize}
%   \item if $t = x$ then $\trans[\dag]{t} = x$, 
%     \begin{math}
%       \substA{t}{\alpha}{v}{u}^\dag
%       = x 
%       = \evalA{\trans[\dag]{t}}{\alpha}{\trans[\dag]{u}}
%     \end{math}
%   \item if $t = \abs{x}{w}$ then 
%     $\trans[\dag]{t} = \abs{x}{w^\dag}$, 
%     \begin{math}
%       \evalA{\trans[\dag]{t}}{\alpha}{\trans[\dag]{u}}
%       = \abs{x}{\evalA{w^\dag}{\alpha}{\trans[\dag]{u}}}
%     \end{math}
%     and 
%     \begin{displaymath}
%       \begin{array}{rcl}
%         \substA{t}{\alpha}{v}{u}^\dag
%         & = & \substA{\abs{x}{w}}{\alpha}{v}{u}^\dag \\
%         & = & \abs{x}{\substA{w}{\alpha}{v}{u}}^\dag \\
%         & = & \abs{x}{\evalA{w^\dag}{\alpha}{\trans[\dag]{u}}} \\
%       \end{array} 
%     \end{displaymath}
  \item if $t = \app{a}{b}$ then 
    \begin{displaymath}
      \begin{array}{rcl}
        \trans{\substA{\app{a}{b}}{\alpha}{v}{u}} 
%         \mathsf{\lm-subst. \ def.} 
        & = & 
        \trans[\dag]{\app{\substA{a}{\alpha}{v}{u}}{\substA{b}{\alpha}{v}{u}}} \\ 
%         \mathsf{def. \ \ref{def:dag-translation}} 
        & = & 
        \mabs{\beta}{\cut{\trans[\dag]{\substA{b}{\alpha}{v}{u}}}{\eabs{y}
            {\cut{\trans[\dag]{\substA{a}{\alpha}{v}{u}}}{\eapp{y}{\beta}}}}} \\ 
%         \mathsf{ind. \ hyp.} 
        & \lreduce{}{0} & 
        \mabs{\beta}{\cut{\evalA{\trans{b}}{\alpha}{\trans[\dag]{u}}}
          {\eabs{y}{\cut{\evalA{\trans{a}}{\alpha}{\trans[\dag]{u}}}{\eapp{y}{\beta}}}}}\\
%         \mathsf{\lmm-subst. \ def.} 
        & = & 
        \evalA{\mabs{\beta}{\cut{\trans{b}}
            {\eabs{y}{\cut{\trans{a}}{\eapp{y}{\beta}}}}}}{\alpha}{\trans{u}} \\ 
%         \mathsf{def. \ \ref{dag-translation}}
        & = & 
        \evalA{\trans[\dag]{\app{a}{b}}}{\alpha}{\trans{u}}
      \end{array}
    \end{displaymath}
%   \item if $t = \mabs{\beta}{c}$ then 
%     $\trans[\dag]{t} = \mabs{\alpha}{\trans[\dag]{c}}$, 
%     \begin{math}
%       \evalA{\trans[\dag]{t}}{\alpha}{\trans[\dag]{u}}
%       = \mabs{\beta}{\evalA{\trans[\dag]{c}}{\alpha}{\trans[\dag]{u}}} 
%     \end{math} 
%     and 
%     \begin{displaymath}
%       \begin{array}{rcl}
%         \substA{t}{\alpha}{v}{u}^\dag 
%         & = & \substA{\mabs{\beta}{c}}{\alpha}{v}{u}^\dag \\ 
%         & = & \mabs{\beta}{\substA{c}{\alpha}{v}{u}}^\dag \\ 
%         & = & \mabs{\beta}{\substA{c}{\alpha}{v}{u}^\dag} \\ 
%         & = & \mabs{\beta}{\substA{\trans[\dag]{c}}{\alpha}{v}{\trans[\dag]{u}}} \\
%       \end{array}
%     \end{displaymath}
  \item if $c = \mapp{\alpha}{w}$ then 
    \begin{displaymath}
      \begin{array}{rcl}
        \trans[\dag]{\substA{\mapp{\alpha}{w}}{\alpha}{v}{u}} 
        % \mathsf{\lm-subst. \ def.} 
        & = & \trans[\dag]{\subst{\mapp{\alpha}{w}}{\alpha}{\happ{\alpha}{u}}} \\
        & = & \trans[\dag]{\happ{\alpha}{u} \hcst{\subst{w}{\alpha}{\happ{\alpha}{u}}}} \\
        & = & \cut{\trans[\dag]{\subst{w}{\alpha}{\happ{\alpha}{u}}}}{\eapp{\trans[\dag]{u}}{\alpha}} \\
        & \lreduce{}{0} & \cut{\subst{\trans[\dag]{w}}{\alpha}{\happ{\alpha}{\trans[\dag]{u}}}}{\eapp{\trans[\dag]{u}}{\alpha}} \\
        & = & \subst{\cut{\trans[\dag]{w}}{\alpha}}{\alpha}{\eapp{\trans[\dag]{u}}{\alpha}} \\
        & = & \evalA{\trans[\dag]{\mapp{\alpha}{w}}}{\alpha}{\trans[\dag]{u}}
      \end{array}
    \end{displaymath}
%   \item if $c = \mapp{\beta}{w}$ then 
%     \begin{displaymath}
%       \begin{array}{rcl}
%         & & \trans[\dag]{\substA{\mapp{\beta}{w}}{\alpha}{v}{u}} \\ 
%         \mathsf{\lm-subst. \ def.} & = & 
%         \trans[\dag]{\mapp{\beta}{\substA{w}{\alpha}{v}{u}}} \\ 
%         \mathsf{def. \ \ref{def:dag-translation}} & = & 
%         \cut{\trans[\dag]{\substA{w}{\alpha}{v}{u}}}{\beta} \\
%         \mathsf{ind. \ hyp.} & \lreduce{}{0} & 
%         \cut{\evalA{\trans[\dag]{w}}{\alpha}{\trans[\dag]{u}}}{\beta} \\ 
%         \mathsf{\lmm-subst. \ def.} & = &
%         \evalA{\cut{\trans{w}}{\beta}}{\alpha}{\trans[\dag]{u}} \\ 
%         \mathsf{def. \ \ref{def:dag-translation}} & = & 
%         \evalA{\trans[\dag]{\mapp{\beta}{w}}}{\alpha}{\trans[\dag]{u}} \\ 
%       \end{array}
%     \end{displaymath}
  \end{itemize}
\end{proof}

\newcommand{\substB}[4]{\subst{#1}{#2}{\hpush[#2]{#4}}}
\newcommand{\evalB}[3]{\subst{#1}{#2}{\eabs{y}{\cut{#3}{\eapp{y}{#2}}}}}

\begin{lemma}
  \label{thm:dag-muprime-subst}
  \begin{math}
    \trans{\substB{t}{\alpha}{v}{u}}
    \lreduce{}{0} 
    \evalB{\trans[\dag]{t}}{\alpha}{y}{\trans[\dag]{u}}
  \end{math}
\end{lemma}

\begin{proof} By induction on $t$.
  \begin{itemize}
%   \item if $t = x$ then $\trans[\dag]{t} = x$, 
%     \begin{math}
%       \substB{\trans[\dag]{t}}{\alpha}{v}{\trans[\dag]{u}} 
%       = x 
%       = \evalB{\trans[\dag]{t}}{\alpha}{y}{\trans[\dag]{u}} 
%     \end{math}
%   \item if $t= \abs{x}{w}$ then 
%     $\trans[\dag]{t} = \abs{x}{w^\dag}$, 
%     \begin{math}
%       \evalB{\trans[\dag]{t}}{\alpha}{y}{\trans[\dag]{u}}
%       = \abs{x}{\evalB{w^\dag}{\alpha}{y}{\trans[\dag]{u}}}
%     \end{math}
%     and 
%     \begin{displaymath}
%       \begin{array}{rcl}
%         \substB{t}{\alpha}{v}{u}^\dag
%         & = & \substB{\abs{x}{w}}{\alpha}{v}{u}^\dag \\
%         & = & \abs{x}{\substB{w}{\alpha}{v}{u}}^\dag \\
%         & = & \abs{x}{\evalB{w^\dag}{\alpha}{y}{\trans[\dag]{u}}} \\
%       \end{array}
%     \end{displaymath}
  \item if $t = \app{a}{b}$ then 
    \begin{displaymath}
      \begin{array}{rcl}
        \trans{\substB{\app{a}{b}}{\alpha}{v}{u}} 
%         \mathsf{\lm-subst. \ def.} 
        & = & 
        \trans[\dag]{\app{\substB{a}{\alpha}{v}{u}}{\substB{b}{\alpha}{v}{u}}} \\ 
%         \mathsf{def. \ \ref{def:dag-translation}} 
        & = & 
        \mabs{\beta}{\cut{\trans[\dag]{\substB{b}{\alpha}{v}{u}}}{\eabs{y}
            {\cut{\trans[\dag]{\substB{a}{\alpha}{v}{u}}}{\eapp{y}{\beta}}}}} \\ 
%         \mathsf{ind. \ hyp.} 
        & \lreduce{}{0} & 
        \mabs{\beta}{\cut{\evalB{\trans{b}}{\alpha}{\trans[\dag]{u}}}
          {\eabs{y}{\cut{\evalB{\trans{a}}{\alpha}{\trans[\dag]{u}}}{\eapp{y}{\beta}}}}}\\
%         \mathsf{\lmm-subst. \ def.} 
        & = & 
        \evalB{\mabs{\beta}{\cut{\trans{b}}
            {\eabs{y}{\cut{\trans{a}}{\eapp{y}{\beta}}}}}}{\alpha}{\trans{u}} \\ 
%         \mathsf{def. \ \ref{dag-translation}}
        & = & 
        \evalB{\trans[\dag]{\app{a}{b}}}{\alpha}{\trans{u}}
      \end{array}
    \end{displaymath}
%   \item if $t = \mabs{\beta}{c}$ then 
%     $\trans[\dag]{t} = \mabs{\alpha}{\trans[\dag]{c}}$, 
%     \begin{math}
%       \evalB{\trans[\dag]{t}}{\alpha}{y}{\trans[\dag]{u}}
%       = \mabs{\beta}{\evalB{\trans[\dag]{c}}{\alpha}{y}{\trans[\dag]{u}}} 
%     \end{math} 
%     and 
%     \begin{displaymath}
%       \begin{array}{rcl}
%         \substB{t}{\alpha}{v}{u}^\dag 
%         & = & \substB{\mabs{\beta}{c}}{\alpha}{v}{u}^\dag \\ 
%         & = & \mabs{\beta}{\substB{c}{\alpha}{v}{u}}^\dag \\ 
%         & = & \mabs{\beta}{\substB{c}{\alpha}{v}{u}^\dag} \\ 
%         & = & \mabs{\beta}{\evalB{\trans[\dag]{c}}{\alpha}{y}{\trans[\dag]{u}}} \\
%       \end{array}
%     \end{displaymath}
  \item if $c = \mapp{\alpha}{w}$ then 
    \begin{displaymath}
      \begin{array}{rcl}
        \trans[\dag]{\substB{\mapp{\alpha}{w}}{\alpha}{v}{u}} 
        % \mathsf{\lm-subst. \ def.} 
        & = & 
        \trans[\dag]{\mapp{\alpha}{\app{\substB{w}{\alpha}{v}{u}}{u}}} \\ 
        % \mathsf{def. \ \ref{def:dag-translation}} 
        & = & 
        \cut{\trans[\dag]{\app{\substB{w}{\alpha}{v}{u}}{u}}}{\alpha} \\ 
        % \mathsf{lem. \ \ref{thm:dag-applicative-sequence}} 
        & \lreduce{}{0} & 
        \cut{\trans[\dag]{\substB{w}{\alpha}{v}{u}}}{\eapp{\trans[\dag]{u}}{\alpha}} \\ 
        % \mathsf{ind. \ hyp.} 
        & \lreduce{}{0} & 
        \cut{\evalB{\trans{w}}{\alpha}{\trans{u}}}{\eapp{\trans{u}}{\alpha}} \\ 
        % \mathsf{\lmm-subst. \ def.} 
        & = & 
        \evalB{\cut{\trans{w}}{\alpha}}{\alpha}{\trans[\dag]{u}} \\ 
        % \mathsf{def. \ \ref{def:dag-translation}} 
        & = & 
        \evalB{\trans[\dag]{\mapp{\alpha}{w}}}{\alpha}{\trans[\dag]{u}}
      \end{array}
    \end{displaymath}
%   \item if $c = \mapp{\beta}{w}$ then 
%     \begin{displaymath}
%       \begin{array}{rcl}
%         \trans[\dag]{\substB{\mapp{\beta}{w}}{\alpha}{v}{u}} 
%         % \mathsf{\lm-subst. \ def.} 
%         & = & 
%         \trans[\dag]{\mapp{\beta}{\substB{w}{\alpha}{v}{u}}} \\ 
%         % \mathsf{def. \ \ref{def:dag-translation}} 
%         & = & 
%         \cut{\trans[\dag]{\substB{w}{\alpha}{v}{u}}}{\beta} \\
%         % \mathsf{ind. \ hyp.} 
%         & \lreduce{}{0} & 
%         \cut{\evalB{\trans[\dag]{w}}{\alpha}{\trans[\dag]{u}}}{\beta} \\ 
%         % \mathsf{\lmm-subst. \ def.} 
%         & = &
%         \evalB{\cut{\trans{w}}{\beta}}{\alpha}{\trans[\dag]{u}} \\ 
%         % \mathsf{def. \ \ref{def:dag-translation}} 
%         & = & 
%         \evalB{\trans[\dag]{\mapp{\beta}{w}}}{\alpha}{\trans[\dag]{u}} \\ 
%       \end{array}
%     \end{displaymath}
  \end{itemize}
\end{proof}

\newcommand{\substC}[1]{\subst{#1}{\alpha}{\beta}}
\newcommand{\evalC}[1]{\subst{#1}{\alpha}{\beta}}

\begin{lemma}
  \label{thm:dag-theta-subst}
  \begin{math}
    \trans{\substC{t}} = \evalC{\trans{t}} 
  \end{math}
\end{lemma}

\begin{proof} By induction on $t$. 
  \begin{itemize}
  \item if $c = \mapp{\alpha}{u}$ then 
    \begin{displaymath}
      \begin{array}{rcl}
        \trans{\substC{\mapp{\alpha}{u}}} 
        % \mathsf{\lm-subst. \ def.} 
        & = & 
        \mapp{\beta}{\trans{\substC{u}}} \\
        % \mathsf{def. \ \ref{def:dag-translation}} 
        & = & 
        \cut{\trans{\substC{u}}}{\beta} \\
        % \mathsf{ind. \ hyp.} 
        & = & 
        \cut{\evalC{\trans{u}}}{\beta} \\ 
        % \mathsf{\lmm-subst. \ def.} 
        & = & 
        \evalC{\cut{\trans{u}}{\alpha}} \\ 
        % \mathsf{def. \ \ref{def:dag-translation}} 
        & = & 
        \evalC{\trans{\mapp{\alpha}{u}}}
      \end{array}
    \end{displaymath}
%   \item if $c = \mapp{\gamma}{u}$ then 
%     \begin{displaymath}
%       \begin{array}{rcl}
%         & & \trans{\substC{\mapp{\gamma}{u}}} \\ 
%         \mathsf{\lm-subst. \ def.} & = & 
%         \trans{\mapp{\gamma}{\substC{u}}} \\
%         \mathsf{def. \ \ref{def:dag-translation}} & = & 
%         \cut{\beta}{\trans{\substC{u}}} \\
%         \mathsf{ind. \ hyp.} & = & 
%         \cut{\evalC{\trans{u}}}{\gamma} \\ 
%         \mathsf{\lmm-subst. \ def.} & = & 
%         \evalC{\cut{\trans{u}}{\gamma}} \\ 
%         \mathsf{def. \ \ref{def:dag-translation}} & = & 
%         \evalC{\trans{\mapp{\gamma}{u}}}
%       \end{array}
%     \end{displaymath}
  \end{itemize}
\end{proof}

\begin{theorem}[simulation of the $\lm$-calculus by the $\lmm$-calculus]
  \label{thm:lm-lmm-simulation}
  \begin{displaymath}
    t \reduce{\gamma}{1} v 
    \ \Longrightarrow \ 
    \exists u \ \trans[\dag]{t} \reduce{}{0} u \ \land \  \trans{v} \lreduce{}{0} u 
\end{displaymath}
\end{theorem}

\begin{proof} By cases on $\gamma$.
  \begin{itemize}
  \item if $\gamma = \beta$ then 
    \begin{displaymath}
      \begin{array}{rcl}
        \trans{\app{\abs{x}{u}}{v}} 
        % \mathsf{def. \ \ref{def:dag-translation}} 
        & = & 
        \mabs{\beta}{\cut{\trans[\dag]{v}}{\eabs{y}
            {\cut{\abs{x}{\trans[\dag]{u}}}{\eapp{y}{\beta}}}}} \\
        & \lreduce{\beta}{1} & 
        \mabs{\beta}{\cut{\trans[\dag]{v}}{\eabs{y}
            {\cut{y}{\eabs{x}{\cut{\trans[\dag]{u}}{\beta}}}}}} \\
        & \lreduce{\mutilde}{1} & 
        \mabs{\beta}{\cut{\trans{v}}{\eabs{x}{\cut{\trans{u}}{\beta}}}} \\
        & \reduce{\mutilde}{1} & 
        \mabs{\beta}{\cut{\subst{\trans[\dag]{u}}{x}{\trans[\dag]{v}}}{\beta}} \\
        & \lreduce{\theta}{1} & \subst{\trans[\dag]{u}}{x}{\trans[\dag]{v}} \\ 
        % \mathsf{lem. \ \ref{thm:dag-beta-subst}} 
        & = & \trans{\subst{u}{x}{v}} 
      \end{array}
    \end{displaymath}
  \item if $\gamma = \mu$ then 
    \begin{displaymath}
      \begin{array}{rcl}
        \trans{\app{\mabs{\alpha}{c}}{v}} 
        % \mathsf{def. \ \ref{def:dag-translation}} 
        & = & 
        \mabs{\beta}{\cut{\trans{v}}{\eabs{y}
            {\cut{\mabs{\alpha}{\trans{c}}}{\eapp{y}{\beta}}}}} \\ 
        & \lreduce{\mutilde}{1} & 
        \mabs{\alpha}{\cut{\mabs{\alpha}{\trans{c}}}{\eapp{\trans{v}}{\alpha}}} \\ 
        & \reduce{\mu}{1} & 
        \mabs{\alpha}{\evalA{\trans{c}}{\alpha}{\trans{v}}} \\ 
        % \mathsf{lem. \ \ref{thm:dag-mu-subst}} 
        & \approx & 
        \trans{\mabs{\alpha}{\substA{c}{\alpha}{u}{v}}}
      \end{array}
    \end{displaymath}
  \item if $\gamma = \mu'$ then 
    \begin{displaymath}
      \begin{array}{rcl}
        \trans{\app{v}{\mabs{\alpha}{c}}} 
        % \mathsf{def. \ \ref{def:dag-translation}} 
        & = & 
        \mabs{\beta}{\cut{\mabs{\alpha}{\trans{c}}}
          {\eabs{y}{\cut{\trans{v}}{\eapp{y}{\beta}}}}} \\
        & \reduce{\mu}{1} & 
        \mabs{\alpha}{\evalB{\trans{c}}{\alpha}{\trans{v}}} \\
        % \mathsf{lem. \ \ref{thm:dag-muprime-subst}} 
        & \approx & 
        \trans{\mabs{\alpha}{\substB{c}{\alpha}{u}{v}}} 
      \end{array}
    \end{displaymath}
  \item if $\gamma = \rho$ then 
    \begin{displaymath}
      \begin{array}{rcl}
        \trans{\mapp{\beta}{\mabs{\alpha}{c}}} 
        % \mathsf{def. \ \ref{def:dag-translation}} 
        & = & 
        \cut{\mabs{\alpha}{\trans{c}}}{\beta} \\ 
        & \lreduce{\mu}{1} & \subst{\trans{c}}{\alpha}{\beta} \\
        % \mathsf{lem. \ \ref{thm:dag-theta-subst}} 
        & = & 
        \trans{\subst{c}{\alpha}{\beta}} 
      \end{array}
    \end{displaymath}
  \item if $\gamma = \theta$ then 
    \begin{math}
      \trans{\mabs{\delta}{\mapp{\delta}{t}}} 
      = \mabs{\delta}{\cut{\trans{t}}{\delta}} \
      \lreduce{\theta}{1} \trans{t} 
    \end{math}
  \end{itemize}
\end{proof}

\begin{corollary}[call-by-name case]
  \label{thm:cbn-lm-lmm-simulation}
  \begin{math}
    t \reduce{n}{1} v 
    \ \Longrightarrow \ 
    \exists u \ \trans[\dag]{t} \reduce{n}{0} u \ \land \  \trans{v} \lreduce{n}{0} u 
  \end{math}
\end{corollary}

\begin{proof} By cases on $\beta$ and $\mu$-rules. 
  
  $\app{\abs{x}{u}}{v}$ is $\beta$-reduced in call-by-name without any 
  restriction. It is simulated in the $\lmm$-calculus by a $\mutilde$-reduction. 
  The latter is in call-by-name without any restriction too. 

  $\app{\mabs{\alpha}{c}}{v}$ is $\mu$-reduced in call-by-name without 
  any restriction. It is simulated in the $\lmm$-calculus 
  by a $\mu$-reduction. The latter is in call-by-name if $\eapp{\trans[\dag]{v}}{\alpha}$ 
  is a stack. It is the case by  definition \ref{def:dag-translation}. 
\end{proof}

\begin{corollary}[call-by-value case]
  \label{thm:cbv-lm-lmm-simulation}
  \begin{math}
    t \reduce{v}{1} v 
    \ \Longrightarrow \ 
    \exists u \ \trans[\dag]{t} \reduce{v}{0} u \ \land \  \trans{v} \lreduce{v}{0} u 
  \end{math}
\end{corollary}

\begin{proof} By cases on $\beta$, $\mu$ and $\mu'$-rules. 

  $\app{\abs{x}{u}}{v}$ is $\beta$-reduced in call-by-value 
  if $v$ is a value. It is simulated in the $\lmm$-calculus by a $\mutilde$-reduction. 
  The latter is in call-by-value if $\trans[\dag]{v}$ is a value. 
  It is the case by the definition of $\lmm_Q$. 

  $\app{\mabs{\alpha}{c}}{v}$ is $\mu$-reduced in call-by-value 
  if $v$ is a value. It is simulated in the $\lmm$-calculus by a $\mu$-reduction. 
  The latter is in call-by-value without any restriction. 

  $\app{v}{\mabs{\alpha}{c}}$ is $\mu'$-reduced in call-by-value 
  without any restriction. 
  It is simulated in the $\lmm$-calculus by a $\mu$-reduction. 
  The latter is in call-by-value without any restriction as well.
\end{proof}

The $\lmm$-simulation by the $\lm$-calculus requires preliminary lemmas showing 
that $\trans[\circ]{(~)}$ commutes over $\lm$ and $\lmm$-substitutions. 
Each proof consists of 
\begin{itemize}
\item expanding the $\lmm$-substitution
\item expanding the definition of $\trans[\circ]{(~)}$
\item applying the induction hypothesis if necessary
\item factorising the $\lm$-substitution
\item factorising the definition of $\trans[\circ]{(~)}$
\end{itemize} 

\begin{lemma}
  \label{thm:circ-beta-subst}
  \begin{math}
    \trans[\circ]{\subst{t}{x}{u}} = \subst{\trans[\circ]{t}}{x}{\trans[\circ]{u}} 
  \end{math}
\end{lemma}

\begin{proof} By induction on $t$.
  \begin{itemize}
  \item if $t = x$ then 
    \begin{math}
      \trans[\circ]{\subst{x}{x}{u}} 
      = \trans[\circ]{u} 
      = \subst{\trans[\circ]{x}}{x}{\trans[\circ]{u}}
    \end{math}
%     \begin{displaymath}
%       \begin{array}{rcl}
%         & & \trans[\circ]{\subst{x}{x}{u}} \\ 
%         \mathsf{\lmm-subst. \ def.} & = & \trans[\circ]{u} \\ 
%         \mathsf{\lm-subst. \ def.} & = & \subst{\trans[\circ]{x}}{x}{\trans[\circ]{u}}
%       \end{array}
%     \end{displaymath}
  \item if $t = y$ then 
    \begin{math}
      \trans[\circ]{\subst{y}{x}{u}}
      = y 
      = \subst{\trans[\circ]{y}}{x}{\trans[\circ]{u}}
    \end{math}
%     \begin{displaymath}
%       \begin{array}{rcl}
%         & & \trans[\circ]{\subst{y}{x}{u}} \\ 
%         \mathsf{\lmm-subst. \ def.} & = & \trans[\circ]{y} \\ 
%         \mathsf{\lm-subst. \ def.} & = & \subst{\trans[\circ]{y}}{x}{\trans[\circ]{u}}
%       \end{array}
%     \end{displaymath}
%   \item if $t = \abs{y}{v}$ then 
%     $\trans[\circ]{t} = \abs{y}{\trans[\circ]{v}}$, 
%     $\subst{\trans[\circ]{t}}{x}{\trans[\circ]{u}} = \abs{y}{\subst{\trans[\circ]{v}}{x}{\trans[\circ]{u}}}$ and 
%     \begin{displaymath}
%       \begin{array}{rcl}
%         \subst{t}{x}{u}^\circ 
%         & = & \subst{\abs{y}{v}}{x}{u}^\circ \\
%         & = & \abs{y}{\subst{v}{x}{u}}^\circ \\
%         & = & \abs{y}{\subst{v}{x}{u}^\circ} \\
%         & = & \abs{y}{\subst{\trans[\circ]{v}}{x}{\trans[\circ]{u}}} \\
%       \end{array}
%     \end{displaymath}
%   \item if $t = \mabs{\alpha}{c}$ then 
%     $\trans[\circ]{t} = \mabs{\alpha}{\trans[\circ]{c}}$, 
%     $\subst{\trans[\circ]{t}}{x}{\trans[\circ]{u}} = \mabs{\alpha}{\subst{\trans[\circ]{c}}{x}{\trans[\circ]{u}}}$ and 
%     \begin{displaymath}
%       \begin{array}{rcl}
%         \subst{t}{x}{u}^\circ 
%         & = & \subst{\mabs{\alpha}{c}}{x}{u}^\circ \\
%         & = & \mabs{\alpha}{\subst{c}{x}{u}}^\circ \\
%         & = & \mabs{\alpha}{\subst{c}{x}{u}^\circ} \\
%         & = & \mabs{\alpha}{\subst{\trans[\circ]{c}}{x}{\trans[\circ]{u}}} \\
%       \end{array}
%     \end{displaymath}
  \item if $t = \cut{t}{e}$ then 
    \begin{displaymath}
      \begin{array}{rcl}
        \trans[\circ]{\subst{\cut{t}{e}}{x}{u}} 
        % \mathsf{\lmm-subst. \ def.} 
        & = & 
        \trans[\circ]{\cut{\subst{t}{x}{u}}{\subst{e}{x}{u}}} \\ 
        % \mathsf{def. \ \ref{def:circ-translation}} 
        & = & 
        \trans[\circ]{\subst{e}{x}{u}} \hcst{\trans[\circ]{\subst{t}{x}{u}}} \\ 
        % \mathsf{ind. \ hyp.} 
        & = & 
        \subst{\trans[\circ]{e}}{x}{\trans[\circ]{u}} 
        \hcst{\subst{\trans[\circ]{t}}{x}{\trans[\circ]{u}}} \\ 
        % \mathsf{\lm-subst. \ def.} 
        & = & 
        \subst{\trans[\circ]{e} \hcst{\trans[\circ]{t}}}{x}{\trans[\circ]{u}} \\ 
        % \mathsf{def. \ \ref{def:circ-translation}} 
        & = & 
        \subst{\trans[\circ]{\cut{t}{e}}}{x}{\trans[\circ]{u}} 
      \end{array}
    \end{displaymath}
  \end{itemize}
\end{proof}

\newcommand{\SubstA}[3]{\subst{#1}{#2}{#3}}

\begin{lemma}
  \label{thm:circ-mu-subst}
  \begin{math}
    \trans[\circ]{\subst{t}{\alpha}{h}} = 
    \SubstA{\trans[\circ]{t}}{\alpha}{\trans[\circ]{h}} 
  \end{math}
\end{lemma}

\begin{proof} By induction on $t$.  
  \begin{itemize}
%   \item if $t = x$ then 
%     $\trans[\circ]{t} = x$ and 
%     $\subst{t}{\alpha}{h}^\circ = x^\circ = \SubstA{\trans[\circ]{t}}{\alpha}{\trans[\circ]{h}}$
%   \item if $t = \abs{x}{u}$ then 
%     $\trans[\circ]{t} = \abs{x}{u ^\circ}$, 
%     $\SubstA{\trans[\circ]{t}}{\alpha}{\trans[\circ]{h}} = \abs{x}{\SubstA{\trans[\circ]{u}}{\alpha}{\trans[\circ]{h}}}$ and 
%     \begin{displaymath}
%       \begin{array}{rcl}
%         \subst{t}{\alpha}{h}^\circ 
%         & = & \subst{\abs{x}{u}}{\alpha}{h}^\circ \\
%         & = & \abs{x}{\subst{u}{\alpha}{h}}^\circ \\
%         & = & \abs{x}{\subst{u}{\alpha}{h}^\circ} \\
%         & = & \abs{x}{\SubstA{\trans[\circ]{u}}{\alpha}{\trans[\circ]{h}}} \\
%       \end{array}
%     \end{displaymath}
%   \item if $t = \mabs{\beta}{c}$ then 
%     $\trans[\circ]{t} = \mabs{\beta}{c ^\circ}$, 
%     $\SubstA{\trans[\circ]{t}}{\alpha}{\trans[\circ]{h}} = \mabs{\beta}{\SubstA{\trans[\circ]{c}}{\alpha}{\trans[\circ]{h}}}$ 
%     and 
%     \begin{displaymath}
%       \begin{array}{rcl}
%         \subst{t}{\alpha}{h}^\circ 
%         & = & \subst{\mabs{\beta}{c}}{\alpha}{h}^\circ \\
%         & = & \mabs{\beta}{\subst{c}{\alpha}{h}}^\circ \\
%         & = & \mabs{\beta}{\subst{c}{\alpha}{h}^\circ} \\
%         & = & \mabs{\beta}{\SubstA{\trans[\circ]{c}}{\alpha}{\trans[\circ]{h}}} \\
%       \end{array}
%     \end{displaymath}
  \item if $c = \cut{t}{e}$ then 
    \begin{displaymath}
      \begin{array}{rcl}
        \trans[\circ]{\subst{\cut{t}{e}}{\alpha}{h}} 
        % \mathsf{\lmm-subst. \ def.} 
        & = & 
        \trans[\circ]{\cut{\subst{t}{\alpha}{h}}{\subst{e}{\alpha}{h}}} \\ 
        % \mathsf{def. \ \ref{def:circ-translation}} 
        & = & 
        \trans[\circ]{\subst{e}{\alpha}{h}} \hcst{\trans[\circ]{\subst{t}{\alpha}{h}}} \\ 
        % \mathsf{ind. \ hyp.} 
        & = & 
        \SubstA{\trans[\circ]{e}}{\alpha}{\trans[\circ]{h}} 
        \hcst{\SubstA{\trans[\circ]{t}}{\alpha}{\trans[\circ]{h}}} \\ 
        % \mathsf{\lm-subst. \ def.} 
        & = & 
        \SubstA{\trans[\circ]{e} \hcst{\trans[\circ]{t}}}{\alpha}{\trans[\circ]{h}} \\ 
        % \mathsf{def. \ \ref{def:circ-translation}} 
        & = & 
        \SubstA{\trans[\circ]{\cut{t}{e}}}{\alpha}{\trans[\circ]{h}}
      \end{array}
    \end{displaymath}
  \item if $e = \evar{\alpha}$ then 
    \begin{math}
      \trans[\circ]{\subst{\evar{\alpha}}{\alpha}{h}} 
      = \trans[\circ]{h} 
      = \SubstA{\trans[\circ]{\alpha}}{\alpha}{\trans[\circ]{h}} 
    \end{math}
%     \begin{displaymath}
%       \begin{array}{rcl}
%         & & \trans[\circ]{\subst{\evar{\alpha}}{\alpha}{h}} \\ 
%         \mathsf{\lmm-subst. \ def.} & = & 
%         \trans[\circ]{h} \\
%         \mathsf{\lm-subst. \ def.} & = & 
%         \SubstA{\hvar{\alpha}}{\alpha}{\trans[\circ]{h}} \\
%         \mathsf{def. \ \ref{def:circ-translation}} & = & 
%         \SubstA{\trans[\circ]{\alpha}}{\alpha}{\trans[\circ]{h}} 
%       \end{array}
%     \end{displaymath}
  \item if $e = \evar{\beta}$ then 
    \begin{math}
      \trans[\circ]{\subst{\evar{\beta}}{\alpha}{h}} 
      = \trans[\circ]{\beta} 
      = \SubstA{\trans[\circ]{\beta}}{\alpha}{\trans[\circ]{h}} 
    \end{math}
 \end{itemize}
\end{proof}

\begin{theorem}[simulation of the $\lmm$-calculus by the $\lm$-calculus]
  \label{thm:lmm-lm-simulation}
  \begin{displaymath}
    t \reduce{\gamma}{1} v 
    \Longrightarrow 
    \exists u \ \trans[\circ]{t} \reduce{}{0} u \lreduce{}{0} \trans[\circ]{v}
  \end{displaymath}
\end{theorem}

\begin{proof} By cases on $\gamma$. 
  \begin{itemize}
  \item if $\gamma = \beta'$ then 
    \begin{displaymath}
      \begin{array}{rcl}
        \trans[\circ]{\cut{\abs{x}{u}}{\eapp{v}{e}}} 
        % \mathsf{def. \ \ref{def:circ-translation}} 
        & = & 
        \hpop{\trans[\circ]{e}}{\trans[\circ]{v}} \hcst{\abs{x}{\trans[\circ]{u}}} \\ 
        % \mathsf{def. \ \ref{def:lm-term-context-cut}} 
        & = & 
        \trans[\circ]{e} \hcst{\app{\abs{x}{\trans[\circ]{u}}}{\trans[\circ]{v}}} \\ 
        & \reduce{\beta}{1} & 
        \trans[\circ]{e} \hcst{\subst{\trans[\circ]{u}}{x}{\trans[\circ]{v}}} \\ 
        % \mathsf{lem. \ \ref{thm:circ-beta-subst}} 
        & = & 
        \trans[\circ]{e} \hcst{\trans[\circ]{\subst{u}{x}{v}}} \\ 
        % \mathsf{def. \ \ref{def:circ-translation}} 
        & = & 
        \trans[\circ]{\cut{\subst{u}{x}{v}}{e}}
      \end{array}
    \end{displaymath}
  \item if $\gamma = \mu$ then 
    \begin{displaymath}
      \begin{array}{rcl}
        \trans[\circ]{\cut{\mabs{\alpha}{c}}{e}} 
        % \mathsf{def. \ \ref{def:circ-translation}} 
        & = & 
        \trans[\circ]{e} \hcst{\mabs{\alpha}{\trans[\circ]{c}}} \\ 
        % \mathsf{lem. \ \ref{thm:lm-cut-subst}} 
        & \reduce{}{0} & 
        \SubstA{\trans[\circ]{c}}{\alpha}{\trans[\circ]{e}} \\ 
        % \mathsf{lem. \ \ref{thm:circ-mu-subst}} 
        & = & 
        \trans[\circ]{\subst{c}{\alpha}{e}} 
      \end{array}
    \end{displaymath}
  \item if $\gamma = \mutilde$ then 
    \begin{displaymath}
      \begin{array}{rcl}
        \trans[\circ]{\cut{t}{\eabs{x}{c}}} 
        % \mathsf{def.  \ \ref{def:circ-translation} \ and \ def. \ \ref{def:lm-term-context-cut}} 
        & = & 
        \mapp{\beta}{\app{\abs{x}{\mabs{\delta}{\trans[\delta]{c}}}}{\trans[\circ]{t}}} \\
        & \reduce{\beta}{1} & 
        \mapp{\beta}{\mabs{\delta}{\subst{\trans[\circ]{c}}{x}{\trans[\circ]{t}}}} \\ 
        & \lreduce{\rho}{1} & 
        \subst{\trans[\circ]{c}}{x}{\trans[\circ]{t}} \\ 
        % \mathsf{lem. \ \ref{thm:circ-beta-subst}} 
        & = & 
        \trans[\circ]{\subst{c}{x}{t}} 
      \end{array}
    \end{displaymath}
  \item if $\gamma = \theta$ then 
    \begin{math}
      \trans[\circ]{\mabs{\delta}{\cut{t}{\delta}}} 
      = 
      \mabs{\delta}{\mapp{\delta}{\trans[\circ]{t}}} \
      \lreduce{\theta}{1} \trans[\circ]{t} 
    \end{math}
  \end{itemize}
\end{proof}

\begin{corollary}[call-by-name case]
  \label{thm:cbn-lmm-lm-simulation}
  \begin{math}
    t \reduce{n}{1} v 
    \Longrightarrow 
    \exists u \ \trans[\circ]{t} \reduce{n}{0} u \lreduce{n}{0} \trans[\circ]{v}
  \end{math}
\end{corollary}

\begin{proof} By cases on $\beta'$, $\mu$ and $\mutilde$-rules. 

  $\cut{\abs{x}{u}}{\eapp{v}{e}}$ is $\beta'$-reduced in call-by-name 
  without any restriction. It is simulated in the $\lm$-calculus 
  by a $\beta$-reduction. The latter is in call-by-name without 
  any restriction too. 

  $\cut{\mabs{\alpha}{c}}{e}$ is $\mu$-reduced in call-by-name 
  if $e \neq \eabs{x}{c'}$ else it were $\mutilde$-reduced. 
  It is simulated in the $\lm$-calculus with the help of lemma 
  \ref{thm:lm-cut-subst}. The latter is in call-by-name 
  if $\trans[\circ]{e} \neq \hpush{t}$ i.e. if $e \neq \eabs{x}{c'}$. 
  It is the case by definition \ref{def:circ-translation}. 

  $\cut{t}{\eabs{x}{c}}$ is $\mutilde$-reduced in call-by-name 
  without any restriction. It is simulated in the $\lm$-calculus 
  by a $\beta$-reduction. The latter is in call-by-name 
  without any restriction as well. 
\end{proof}

\begin{corollary}[call-by-value case]
  \label{thm:cbv-lmm-lm-simulation}
  \begin{math}
    t \reduce{v}{1} v 
    \Longrightarrow 
    \exists u \ \trans[\circ]{t} \reduce{v}{0} u \lreduce{v}{0} \trans[\circ]{v}
  \end{math}
\end{corollary}

\begin{proof} By cases on $\beta'$, $\mu$ and $\mutilde$-rules. 

  $\cut{\abs{x}{u}}{\eapp{v}{e}}$ is $\beta'$-reduced in call-by-value 
  if $v$ is a value. It is simulated in the $\lm$-calculus 
  by a $\beta$-reduction. The latter is in call-by-value 
  if $\trans[\circ]{v}$ is a value. 
  It is the case by definition \ref{def:circ-translation}. 

  $\cut{\mabs{\alpha}{c}}{e}$ is $\mu$-reduced in call-by-value 
  if $e$ is either a $\mu$-variable or a context of the form $\eapp{v}{h}$ 
  where $v$ is a value or a $\mu$-abstraction by the definition of $\lmm_Q$. 
  It is simulated in the $\lm$-calculus with the help of lemma 
  \ref{thm:lm-cut-subst}. The latter is in call-by-value 
  if $\trans[\circ]{v}$ is a value in a context 
  of the form $\happ{\trans[\circ]{v}}{\trans[\circ]{h}}$ 
  i.e. if $v$ is a value in a $\eapp{v}{h}$ context. 
  It is the case by definition \ref{def:circ-translation}. 

  $\cut{t}{\eabs{x}{c}}$ is $\mutilde$-reduced in call-by-value 
  if $t$ is a value. It is simulated in the $\lm$-calculus 
  by a $\beta$-reduction. The latter is in call-by-value 
  if $\trans[\circ]{t}$ is a value. It is the case by definition 
  \ref{def:circ-translation}. 
\end{proof}

\section{Conclusion}
\label{sec:conclusion}

Analysis of the $\lm$ and $\lmm$-calculi has shown their computational equivalence. 
It holds for undirected evaluations of pure calculi 
(see theorems \ref{thm:lm-lmm-simulation} and \ref{thm:lmm-lm-simulation}). 
This result is then easily obtained for call-by-name and call-by-value evaluations 
(see corollaries \ref{thm:cbn-lm-lmm-simulation}, \ref{thm:cbv-lm-lmm-simulation}, 
\ref{thm:cbn-lmm-lm-simulation} and \ref{thm:cbv-lmm-lm-simulation}). 
It concerns the simple type system too 
(see lemmas \ref{thm:dag-type-compatible} and \ref{thm:circ-type-compatible}). 

The simulation of the $\lmm$-calculus by the $\lm$-calculus 
is smoother than the simulation of the $\lm$-calculus by the $\lmm$-calculus. 
The first is obtained with the help of linear \emph{reductions} whereas 
the second is obtained with the help of linear \emph{expansions}. 
% The following figure 
% illustrates theorems \ref{thm:lm-lmm-simulation} and \ref{thm:lmm-lm-simulation}:

% \centerline{
%   \hfill
%   \xygraph
%   {
%     []{t}="1"
%     -@{>}[r(2.0)]{v}="2"_{}^{}
%     -@{>}[d]{\trans[\dag]{v}}="3"_{}^{}
%     "3" "2" "1"
%     -@{>}[d]{\trans[\dag]{t}}="4"_{}^{}
%     -@{>}[r]{u}="5"_{}^{}
%     -@{(~)}"3"_{}^{}
%     -@{~>}"5"_{}^{}
%   }
%   \hfill
%   \xygraph
%   {
%     []{t}="10"
%     -@{>}[r(2.0)]{v}="20"_{}^{}
%     -@{>}[d]{\trans[\circ]{v}}="30"_{}^{}
%     "30" "20" "10"
%     -@{>}[d]{\trans[\circ]{t}}="40"_{}^{}
%     -@{>}[r]{u}="50"_{}^{}
%     -@{~>}"30"_{}^{}
%   }
%   \hfill
% }

This work can be extended in three ways. The first consists of proving 
the same results for the call-by-value evaluation of the $\lm$-calculus 
defined in \cite{OngStewart97}. The second consists of defining CPS translations 
to $\lambda$-calculus in order to complete \cite{CurienHerbelin00}. 
The third consists of extending the type system to the other logical constants.

\end{document}